\documentclass[english]{article}
\usepackage[T1]{fontenc}
\usepackage[latin1]{inputenc}
\usepackage{babel}
\usepackage{graphics}
\usepackage{psfrag}
 
 
\newtheorem{prop}{Proposition}[section]
\newtheorem{lemma}[prop]{Lemma}
\newtheorem{theorem}[prop]{Theorem}
\newtheorem{defi}[prop]{Definition}
\newtheorem{rem}[prop]{\em Remark\/}
\newtheorem{cor}[prop]{Corollary}
\newtheorem{cond}[prop]{Conditions}
\newtheorem{ex}[prop]{\em Example\/}
 
\newcommand{\bp}{\begin{prop}}
\newcommand{\bl}{\begin{lemma}}
\newcommand{\bt}{\begin{theorem}}
\newcommand{\bd}{\begin{defi}\rm}
\newcommand{\br}{\begin{rem}\rm}
\newcommand{\be}{\begin{equation}}
\newcommand{\bea}{\begin{eqnarray}}
\newcommand{\bpr}{\begin{proof}}
\newcommand{\bc}{\begin{cor}}
\newcommand{\bco}{\begin{cond}}
\newcommand{\bex}{\begin{ex}}

\newcommand{\ep}{\end{prop}}
\newcommand{\el}{\end{lemma}}
\newcommand{\et}{\end{theorem}}
\newcommand{\ed}{\end{defi}}
\newcommand{\er}{\end{rem}}
\newcommand{\ee}{\end{equation}}
\newcommand{\eea}{\end{eqnarray}}
\newcommand{\epr}{\end{proof}}
\newcommand{\ec}{\end{cor}}
\newcommand{\eco}{\end{cond}}
\newcommand{\eex}{\end{ex}}

\newcommand{\nn}{ \nonumber \\ }
\newcommand{\pr}{{\em Proof:\ }}
\newcommand{\qed}{\vrule height 5pt width 5pt depth 0pt}

\newcommand{\sla}{s_L^{_A}}
\newcommand{\tla}{t_L^{_A}}
\newcommand{\sra}{s_R^{_A}}
\newcommand{\tra}{t_R^{_A}}
\newcommand{\slb}{s_L^{_B}}
\newcommand{\tlb}{t_L^{_B}}
\newcommand{\srb}{s_R^{_B}}
\newcommand{\trb}{t_R^{_B}}
\newcommand{\gla}{\gamma_L^{_A}}
\newcommand{\gra}{\gamma_R^{_A}}
\newcommand{\glb}{\gamma_L^{_B}}
\newcommand{\grb}{\gamma_R^{_B}}
\newcommand{\pla}{\pi_L^{_A}}
\newcommand{\pra}{\pi_R^{_A}}
\newcommand{\plb}{\pi_L^{_B}}
\newcommand{\prb}{\pi_R^{_B}}

\newcommand{\id}{{\rm id}}
 
\newcommand{\ui}{^{(1)}}
\newcommand{\uii}{^{(2)}}
\newcommand{\di}{_{(1)}}
\newcommand{\dii}{{}_{(2)}}
\newcommand{\ot}{\otimes}
\newcommand{\asr}{{A^R}}
\newcommand{\atr}{{^R\! A}}
\newcommand{\asl}{{_L\! A}}
\newcommand{\atl}{{A_L}}
\newcommand{\asrd}{{{\cal A}^*}}
\newcommand{\atrd}{{^*\!\! {\cal A}}}
\newcommand{\atld}{{{\cal A}_*}}
\newcommand{\asld}{{_*{\cal A}}}
\newcommand{\ci}{\circ}
\newcommand{\ila}{{\cal I}^L({\cal A})}

\renewcommand{\c}{{\cal C}}

\newcommand{\err}{i}
\newcommand{\fsr}{{\ell_R}}
\newcommand{\ftr}{{_R\ell}}

\newcommand{\lu}{\leftharpoonup}
\newcommand{\ld}{\leftharpoondown}
\newcommand{\ru}{\rightharpoonup}
\newcommand{\rd}{\rightharpoondown}
\newcommand{\lsr}{\lambda\!^*}

\newcommand{\phisl}{_* \!\phi}
\newcommand{\phisr}{\phi\!^*}
\newcommand{\phitr}{^*\!\! \phi}
\newcommand{\phitl}{\phi_*}
\renewcommand{\i}{\iota}
\newcommand{\ib}{{\bar{\iota}}}
\newcommand{\x}{\times}
\newcommand{\f}{{\cal F}}
\newcommand{\fdot}{\dot{\cal F}}
\newcommand{\stac}[1]{\stackrel{\ot}{_{_{#1}}}}
\newcommand{\ix}{{\cal X}}
\newcommand{\ixb}{\overline{{\cal X}}}
\newcommand{\re}{\Upsilon}
\renewcommand{\a}{{\rm\bf a}}
\renewcommand{\l}{{\rm\bf l}}
\renewcommand{\r}{{\rm\bf r}}
\newcommand{\ev}{{\rm ev}}
\newcommand{\coev}{{\rm coev}}
\newcommand{\BIM}{{\rm\bf BIM}}

\newcommand{\setc}[1]{\setcounter{equation}{#1}}

\newcommand{\lb}{\label}

\begin{document}

\large
\title{\bf Hopf Algebroid Symmetry of Abstract Frobenius
Extensions of Depth 2}
 
\author{\sc Gabriella B\"ohm $^1$ and  Korn\'el Szlach\'anyi $^2$ \\}

\date{}
 
\maketitle
\normalsize 
\footnotetext[1]{
Research Institute for Particle and Nuclear Physics, Budapest,
H-1525 Budapest 114, P.O.B. 49, Hungary\\
E-mail: BGABR@rmki.kfki.hu\\
  Supported by the Hungarian Scientific Research Fund, OTKA --
  T 034 512, FKFP -- 0043/2001 and the Bolyai J\'anos Fellowship}

\footnotetext[2]{
Research Institute for Particle and Nuclear Physics, Budapest,
H-1525 Budapest 114, P.O.B. 49, Hungary\\
E-mail: SZLACH@rmki.kfki.hu\\
  Supported by the Hungarian Scientific Research Fund,
  OTKA -- T 034 512.}
\vskip 2truecm

\vspace{2cm}

\begin{abstract}  
We study 
Frobenius 1-morphisms $\i$ in an additive bicategory $\c$
satisfying the depth 2 condition. We show that the 2-endomorphism rings
$\c^2(\i\x\ib,\i\x\ib)$ and $\c^2(\ib\x\i,\ib\x\i)$ can be equipped
with dual Hopf algebroid structures. We prove also that a Hopf
algebroid appears as the solution of the above abstract symmetry
problem if and only if it possesses a two sided non-degenerate integral.
\end{abstract}

\bigskip
\bigskip
\bigskip
\bigskip


\section{Introduction}

In the classification of extensions of von Neumann algebras `quantum
groups' and `quantum groupoids' play important role.

The starting point was the result of \cite{Szy,Lo} proving that an
irreducible  finite
index {\em depth 2} (or D2 for short) extension of von Neumann factors
can be realized as a crossed product with a finite dimensional
$C^*$-Hopf algebra -- that is the symmetry of the extension is
described by a {\em finite dimensional $C^*$-Hopf algebra}. About
related results see also 
\cite{Oc,Da,EN,Ya,NW}. The D2 property means that the derived
tower of relative commutants 
$N^{\prime}\cap M\subset N^{\prime}\cap M_1\subset N^{\prime}\cap
M_2\subset  \dots$ of the Jones tower $N\subset M\subset M_1\subset
M_2\dots $ is also a Jones tower \cite{HeOc}.

Allowing for {\em reducible} finite index D2 extensions of $II_1$ von
Neumann factors in \cite{NV1} and of von Neumann algebras with finite
centers in \cite{NSzW} the symmetry of the extension was shown to be
described by a {\em finite dimensional $C^*$-weak Hopf algebra}
introduced in \cite{BSz,N,BNSz}. A Galois correspondence has been
established in \cite{NV2} in the case of finite index finite depth
extensions of $II_1$ factors. The infinite index D2 case has been
treated in \cite{EV}  for arbitrary von Neumann algebras endowed with
a regular operator valued weight.

In \cite{KSz} the analysis was extended to extensions of {\em
rings} $N\to M$. The D2 condition has been generalized as

i) the $N$-$M$ bimodule ${_N M _N }\stac{N} {_N M_ M}$ is a direct
summand in a finite direct sum of copies of ${_N M _M}$,

ii) the $M$-$N$ bimodule ${_M M _N} \stac{N} {_N M _N}$  is a direct
summand in a finite direct sum of copies of ${_M M _N}$.

\noindent
A D2 extension of von Neumann algebras is D2 in this sense.

The symmetry of a D2 ring extension $N\to M$ was shown in \cite{KSz}
to be described by a {\em finite bialgebroid}. That is the
step 2 centralizer rings $A\colon = {\rm End} {_N M _N}$ -- the
endomorphism ring of the bimodule ${_N M _N }$ -- and $B\colon
= (M\stac{N} M)^N$ -- the center the bimodule ${_N M
_N}\stac{N} {_N M _N}$ with multiplication $(b_1\ot b_2)(b^{\prime}_1
\ot b^{\prime}_2)\colon = b^{\prime}_1 b_1\ot b_2 b^{\prime}_2$ --
carry dual bialgebroid structures over the basis $C_N(M)$ -- the
centralizer of $N$ in $M$. The bialgebroid $A$
has a natural left action on $M$. Under the additional
assumption that the module $M_N$ is balanced, the $A$-invariant
subalgebra of $M$ was shown to be $N$. Also $B$ acts naturally on ${\rm
End} {_N M}$ with the invariant subalgebra isomorphic to $M$.

The bialgebroid -- or Takeuchi $\x_L$- 
bialgebra -- is a generalization of the bialgebra in the sense that 
instead of being an $L$-algebra it is a bimodule over a
{\em non-commutative} base ring 
$L$. The bialgebroid definitions \cite{T,Lu,Xu,Sch1,Sz} were shown to be
equivalent in 
\cite{BreMi}. The bialgebroid $A$ over the base $L$ is {\em finite} if the
$L$-module structures on $A$ are finitely generated projective. By the
results in \cite{KSz} in this case the $L$-duals of $A$ also carry
bialgebroid structures. We summarize the basic definitions of the theory of
bialgebroids in Section \ref{preli}. 

In the Section 3 of \cite{HGD} the D2 ring extension $N\to M$ was
supposed to be also {\em Frobenius }. This means the existence of
a Frobenius map $\Phi:{_NM_N}\to N$ possessing a {\em
quasi-basis } $x_j,y_j\in M$ satisfying 
$\sum_j y_j \Phi(x_j m)=m=\sum_j\Phi (m y_j) x_j$ for all $m\in M$. Under
this assumption the above rings
$A$ and $B$ carry more structure then just being
bialgebroids. Namely, fixing a Frobenius map $M\to N$ the bialgebroids
$A$ and $B$ can be equipped 
with antipodes, making them dual {\em Hopf algebroids}.

It is important to emphasize that we use the term {\em Hopf algebroid} in
the sense of \cite{HGD} and not in the sense of \cite{Lu}. Though -- roughly
speaking --  both definitions can be summarized as a `bialgebroid 
with antipode', the two definitions are {\em not} equivalent -- as
it is illustrated by an example in \cite{HGD}. In
\cite{HGD} we give some equivalent definitions of the Hopf
algebroid. These are all common in the feature  that the axioms
are formulated without reference to a particular section of the
canonical projection $A\stac{Z} A\to A\stac{L} A$, which has played
important role in the definition of \cite{Lu}. The most appropriate
form of the definition for the purposes of this paper involves two,
left and right 
bialgebroid structures on the same ring that are connected by the antipode.
In Section \ref{preli} we cite the definition and some
properties of the Hopf algebroid from \cite{HGD}.

A further step of generalization of the `symmetry reconstruction
problem' was done in \cite{Sz} where the symmetry of {\em abstract} D2
extensions was studied. A 1-morphism $\i$ in an additive bicategory $\c$
is an {\em abstract extension} if it possesses a left dual
$\ib$. In the original paper \cite{Sz} the D2 property of $\i$ was
introduced as 

i) the 1-morphism $(\i\x\ib)\x \i$ is a direct
summand in a finite direct sum of copies of $\i$.

\noindent
In later publications (see \cite{KSz} for example) this property has
been renamed as {\em left} D2 property. The condition that

ii) the 1-morphism $(\ib\x\i)\x \ib$ is a direct
summand in a finite direct sum of copies of $\ib$

\noindent
was called {\em right} D2 property. The 1-morphism $\i$ is D2 in this
later sense if both the left and right D2 conditions hold true. In
this paper we use this latter terminology.

The name `abstract extension' is motivated by the basic example:
Observe that for a ring extension $N\to M$ the
bimodule  ${_N M _M}$ is a 1-morphism in the bicategory 
of bimodules possessing a left dual ${_M M _N}$. 
If the extension $N\to M$ is also D2 then the corresponding
1-morphism ${_N M _M}$ satisfies the D2
condition. Hence a D2 ring extension determines a D2 1-morphism in an
additive bicategory -- called an abstract D2 extension.

Motivated by this example in \cite{Sz} the 1-morphism $\i$  in 
$\c$ possessing a left dual $\ib$ is assumed to satisfy the left D2
condition. Then (denoting the horizontal composition in $\c$ by $\x$) the 
ring $A\colon = \c^2(\i\x\ib,\i\x\ib)$ of 2-endomorphisms is shown to
carry a left bialgebroid structure. 
(The right D2 property of $\i$ implies the existence of a right
bialgebroid structure on  $B\colon = \c^2(\ib\x\i,\ib\x\i)$.)
The case when the finite index
property and irreducibility of $\i$ is required and $\c$ is semisimple
and $k$-linear -- leading to semisimple and cosemisimple Hopf algebra
symmetry -- is spelled out 
in detail in \cite{Mu}.


In this paper we generalize the analysis of the Section 3 of
\cite{HGD} to abstract extensions. That is we suppose that $\i$ is a
D2 {\em Frobenius} 1-morphism in an additive bicategory $\c$ with
two sided dual $\ib$. Under this assumption we show that the rings
$A= \c^2(\i\x\ib,\i\x\ib)$  and $B=\c^2(\ib\x\i,\ib\x\i)$ carry dual
Hopf algebroid structures with non-degenerate two sided integrals in
the sense of \cite{HGD}. We also prove that a Hopf algebroid appears
as a solution of this abstract symmetry problem if and only if it
possesses a two sided non-degenerate integral.

It is discussed in \cite{HGD,B} that -- in contrast to (weak) Hopf
algebras 
-- the Hopf algebroid structure on a given left bialgebroid is not
unique. The choice of the antipode maps in Subsection \ref{harm} is a
canonical construction fitting to the Frobenius structure of
$\i$. However, 
this canonical construction is responsible for the existence of the
two sided non-degenerate integral.

This result raises the question how typical it is that we can find 
two sided non-degenerate integrals in Hopf algebroids. As a partial answer,
it is proven in Proposition 5.13 in \cite{HGD} that for any Hopf
algebroid ${\cal A}$ possessing a non-degenerate left integral $\ell$
there exists another Hopf algebroid which is isomorphic to  ${\cal A}$
as a left bialgebroid and in which  $\ell$ is a two sided
non-degenerate integral. In particular, it follows by combining the
Proposition 5.13 in \cite{HGD} with the Theorem 4.2 in \cite{B} that
for any finite dimensional (weak)
Hopf algebra there exists a twisted antipode \cite{ConnMo} leading to
a Hopf algebroid with two sided non-degenerate integral.
It is not known, however, whether any finite Hopf algebroid possesses a
non-degenerate left integral, i.e. no generalization of the
Larson-Sweedler theorem on bialgebroids is known.


The paper is organized as follows: In the Section \ref{preli} we
shortly review some results from \cite{Sw,Sz,KSz} on bialgebroids and
from \cite{HGD} on Hopf algebroids that we are going to use later on. 

In the Subsection \ref{harm} we
analyse the structures of the rings $A=\c^2(\i\x\ib,\i\x\ib)$ and
$B=\c^2(\ib\x\i,\ib\x\i)$ if only the Frobenius property of the
1-morphism $\i$ in the additive bicategory $\c$ is
assumed. Analogously to Section 3 of \cite{HGD} already in this
general situation -- i.e. {\em without} assuming the D2 property --
one can do `Fourier analysis'. That is convolution products can be
introduced both on $A$ and $B$ together with `Fourier transformations'
$A\to B$ relating the original and the convolution products. At this
level of generality one can introduce the maps that are going to be
the antipodes in the D2 case.

In the Subsection \ref{depth2} we impose also the D2 condition on the
1-morphism $\i$. Applying the results of \cite{Sz} this implies the
bialgebroid structures on $A$ and $B$. We prove that these bialgebroid
structures together with the antipodes of Subsection \ref{harm} amount
to dual Hopf algebroids in the sense of \cite{HGD} both of them
possessing two sided non-degenerate integrals.

In the Subsection \ref{inv} we address the question {\em what} Hopf
algebroids appear as symmetries of absratct D2 Frobenius
extensions. Being given a Hopf algebroid ${\cal A}$ possessing a two
sided non-degenerate integral, we construct a D2 Frobenius 1-morphism
${\cal X}$ in the additive bicategory of internal bimodules in its
module category. We show that the symmetry of ${\cal X}$ is described
by ${\cal A}$. Combining this result with the one in Subsection \ref{depth2}
we conclude that the existence of a {\em
two sided non-degenerate integral} \cite{HGD} in the Hopf algebroid is a
sufficient and necessary condition for it to appear as the symmetry of an
abstract D2 Frobenius extension.

In the Appendix we sketch briefly the construction of the bicategory
of internal bimodules in a monoidal category -- used as a tool in
Subsection \ref{inv}.

All rings appearing in the paper are associative and unital. The
category of left/right/bi- modules over the ring $R$ is denoted by
${_R{\cal M}}$/${{\cal M}_R}$/${_R{\cal M}_R}$.

\section{Preliminaries: Bialgebroid and Hopf algebroid}
\lb{preli}
\setc{0}


In this section we summarize our notations and the basic definitions of
bialgebroids and Hopf algebroids that will be used later on. For more
about bialgebroids we refer to the literature
\cite{T,Sch,BreMi,KSz,Sch3,Sz,Sz2} and about Hopf algebroids to
\cite{HGD,B}.

The {\em bialgebroid} \cite{Lu,Xu,Sz} or 
Takeuchi $\x_L$-bialgebra \cite{T}  is a generalization of the
bialgebra to the case of  a non-commutative base ring ring $L$:

\bd A {\em left bialgebroid} ${\cal
A}_L$ consists of the data $(A,L,s_L,t_L,\gamma_L,\pi_L)$. The
$A$ and 
$L$ are associative unital rings, the total and base rings,
respectively. The $s_L:L\to A$ and $t_L:L^{op}\to 
A$ are ring homomorphisms such that the images of $L$ 
under the source map $s_L$ and target map $t_L$ in $A$ commute
making $A$ an $L$-$L$ bimodule:
\be  l\cdot a\cdot l^{\prime}\colon = s_L(l) t_L(l^{\prime}) a .
\lb{elbim} \ee
The bimodule (\ref{elbim}) is denoted by ${_L A_L}$. The triple
$({_L A_L},\gamma_L,\pi_L)$ is a comonoid in ${_L {\cal M}_L}$. Introducing
the Sweedler's notation $\gamma_L(a)\equiv a\di \ot a\dii\in \atl \ot
\asl$ the identities
\bea  a\di t_L(l) \ot a\dii &=& a\di \ot a\dii s_L(l) \lb{cros}\\
      \gamma_L(1_A)&=& 1_A\ot 1_A \\
      \gamma_L(ab)&=&\gamma_L(a) \gamma_L(b) \lb{gmp} \\
      \pi_L(1_A) &=& 1_L \\
      \pi_L\left(a s_L\ci \pi_L(b)\right)=&\pi_L&(ab)=
       \pi_L\left(a t_L\ci \pi_L(b)\right)
\eea
are  required for all $l\in L$ and $a,b\in A$. The requirement 
(\ref{gmp}) makes sense in view of (\ref{cros}).
\ed
If ${\cal A}_L= (A,L,s_L,t_L,\gamma_L,\pi_L)$ is a left bialgebroid
then so is the co-opposite ${\cal A}_{L\ cop}=
(A,L^{op},t_L,s_L,\gamma_L^{op},$ $\pi_L)$ where $\gamma_L^{op}$ denotes
the opposite coproduct $a\mapsto a\dii\ot a\di$. The opposite ${\cal
A}_L^{op}= (A^{op},L,t_L,s_L,\gamma_L,\pi_L)$ is a  right bialgebroid
in the sense of \cite{KSz}: 
\bd A {\em right bialgebroid} ${\cal A}_R$ consists of the data
$(A,R,s_R,t_R,\gamma_R,\pi_R)$. The $A$ and 
$R$ are associative unital rings, the total and base rings,
respectively. The $s_R:R\to A$ and $t_R:R^{op}\to 
A$ are ring homomorphisms such that the images of $R$ 
under the source map $s_R$ and target map $t_R$ in $A$ commute
making $A$ an $R$-$R$ bimodule:
\be  r\cdot a\cdot r^{\prime}\colon = a s_R(r^{\prime}) t_R(r).
\lb{erbim}\ee
The bimodule (\ref{erbim}) is denoted by ${^R\! A^R}$. The triple
$({^R\! A^R},\gamma_R,\pi_R)$ is a comonoid in ${_R {\cal
M}_R}$. Introducing 
the Sweedler's notation $\gamma_R(a)\equiv a\ui \ot a\uii\in \asr \ot
\atr$ the identities
\bea  s_R(r)a\ui  \ot a\uii &=& a\ui \ot  t_R(r) a\uii \nn
      \gamma_R(1_A)&=& 1_A\ot 1_A \nn
      \gamma_R(ab)&=&\gamma_R(a) \gamma_R(b) \nn
      \pi_R(1_A) &=& 1_R \nn
      \pi_R\left(s_R\ci \pi_R(a) b \right)=&\pi_R&(ab)=
       \pi_R\left(t_R\ci \pi_R(a) b \right)
\nonumber\eea
are  required for all $r\in R$ and $a,b\in A$.
\ed
The $L$-actions of the bimodule ${_LA_L}$ are given by left
multiplication and the $R$-actions of the bimodule ${^RA^R}$ are given
by right multiplication. We can define further bimodules:
\bea {^LA^L}:&\qquad& l\cdot a\cdot l^{\prime}\colon = a
t_L(l)s_L(l^{\prime})\lb{lbimr}\\
{_RA_R}:&\qquad& r\cdot a\cdot r^{\prime}\colon = 
s_R(r)t_R(r^{\prime})a.\lb{rbiml}
\eea

In the case of left/right bialgebroids the category $(_A{\cal
M},L,\l,\r,\a)$/$({\cal M}_A,R,\l,\r,\a)$ of left/right $A$-modules
has the  monoidal structure inherited from the bimodule category
${_L{\cal M}_L}$/${_R{\cal M}_R}$. For left/right $A$-modules $M$ and
$N$ the $A$-module structure on $M\ot_L N$/$M\ot_R 
N$ is given by $a\cdot(m\ot n)\colon = a\di \cdot m\ot a\dii\cdot
n$/$(m\ot n)\cdot a\colon = m\cdot a\ui\ot n\cdot a\uii$ -- for $a\in
A$, $m\in M$ and $n\in N$. The monoidal unit  is the left $A$-module
on $L$: $a\cdot l\colon = \pi_L(as_L(l))$/ the right $A$-module on
$R$: $r\cdot a\colon = \pi_R(s_R(r)a)$ -- for $a\in A$, $l\in L$ and
$r\in R$. The forgetful functor ${\Phi}_L:{_A{\cal M}}\to
{_L{\cal M}_L}$/${\Phi}_R:{{\cal M}_A}\to {_R{\cal M}_R}$ is
strong monoidal \cite{Sz}.

The homomorphisms of bialgebroids we use in this paper are the {\em
bialgebroid maps} of \cite{Sz2}:
\bd A {\em left bialgebroid homomorphism} ${\cal A}_L \to {\cal
A}^{\prime}_L$ is a pair of ring homomorphisms $(\Phi:A\to A^{\prime},
\phi: L\to L^{\prime})$ such that 
\bea s^{\prime}_L\ci \phi &=& \Phi\ci s_L\nn
     t^{\prime}_L\ci \phi &=& \Phi\ci t_L\nn
     \pi^{\prime}_L\ci \Phi &=& \phi\ci \pi_L\nn
     \gamma_L^{\prime}\ci \Phi&=& (\Phi\ot \Phi)\ci \gamma_L.
\nonumber\eea
The last condition makes sense since by the first two
conditions $\Phi\ot \Phi$ is a well defined map $\atl \ot \asl\to
{A_{L^{\prime}}}\ot {_{L^{\prime}}}A$.
The pair $(\Phi,\phi)$ is an {\em isomorphism of left bialgebroids} if it is
a left bialgebroid homomorphism such that both $\Phi$ and $\phi$ are 
bijective. 

A right bialgebroid homomorphism (isomorphism) ${\cal A}_R\to {\cal
A}^{\prime}_R$ is a left bialgebroid homomorphism (isomorphism) $({\cal
A}_R)^{op}\to ({\cal A}^{\prime}_R)^{op}$.
\ed


Let ${\cal A}_L$ be a left bialgebroid. The equation (\ref{elbim})
describes two $L$-modules $\atl$ and $\asl$. Their $L$-duals
are the additive groups of $L$-module maps:
\[   \atld\colon = \{ \phitl: \atl \to {L_ L} \} \quad {\rm and} \quad
     \asld\colon = \{ \phisl: \asl \to {_L L} \} \]
where ${_L L}$ stands for the left regular and $L_L$ for the right
regular $L$-module. Both $\atld$ and $\asld$ carry left $A$ module
structures via the transpose of the right regular action of $A$. For
$\phitl\in \atld, \phisl\in\asld$ and $a,b\in A$ we have:
\[   \left(a\ru \phitl\right)(b) =\phitl(ba)\quad {\rm and}\quad
     \left(a\rd \phisl\right)(b) =\phisl(ba). \]
Similarly, in the case of a right bialgebroid ${\cal A}_R$ --  denoting
the left and right regular $R$-modules by $^R R$ 
and $R^R$, respectively, -- the two $R$-dual additive groups
\[  \asrd\colon = \{ \phisr: \asr \to {R^R} \} \quad {\rm and}\quad
    \atrd\colon = \{ \phitr: \atr \to {^R R} \} \]
carry right $A$-module structures:
\[ \left(\phisr \lu a \right)(b) =\phisr(ab) \quad {\rm and}\quad
   \left(\phitr \ld a \right)(b) =\phitr(ab). \]
The comonoid structures can be transposed to give monoid (i.e. ring)
structures to the duals. In the case of a left bialgebroid ${\cal A}_L$ 
\be  \left(\phitl {\psi_*} \right)(a)=
{\psi_*}\left( s_L \ci \phitl (a\di) a\dii \right)\quad {\rm and}\quad
\left(\phisl {_* \psi} \right)(a)=
{_* \psi}\left( t_L \ci \phisl (a\dii) a\di \right) \lb{ldual}\ee
for $\phisl, {_*\psi}\in \asld$, $\phitl,{\psi_*}\in \atld$ and $a\in
A$. 

Similarly, in the case of a right bialgebroid ${\cal A}_R$ 
\be \left(\phisr {\psi^*} \right)(a)=
\phisr\left( a\uii t_R \ci {\psi^*} (a\ui)\right)\quad {\rm and}\quad
 \left(\phitr {^*\! \psi} \right)(a)=
\phitr \left(a\ui  s_R \ci {^* \!\psi} (a\uii)\right)\lb{rdual}\ee
for $\phisr, {\psi^*}\in \asrd$, $\phitr,{^*\psi}\in \atrd$ and $a\in
A$.

In the case of a left bialgebroid ${\cal A}_L$  also the ring $A$ has 
right $\atld$- and   right $\asld$- module sructures:
\[ a\lu \phitl = s_L\ci \phitl (a\di)a\dii \quad {\rm and}\quad
     a\ld \phisl = t_L\ci \phisl (a\dii)a\di \]
for $\phitl\in \atld$, $\phisl\in \asld$ and $a\in A$.

Similarly, in the case of a right bialgebroid ${\cal A}_R$  the ring
$A$ has left $\asrd$- and left $\atrd$ structures:
\[ \phisr\ru a = a\uii t_R\ci \phisr(a\ui )\quad {\rm and}\quad
     \phitr\rd a = a\ui  s_R\ci \phitr(a\uii) \]
for $\phisr\in \asrd$,  $\phitr\in \atrd$ and $a\in A$.

As it was proven in \cite{KSz}, if the $L$ ($R$) module structure
on $A$ is finitely generated projective then the corresponding dual
has also a bialgebroid structure.

\bigskip


The total ring of a Hopf algebroid carries eight canonical module 
structures over the base ring -- modules of the kind (\ref{elbim}),
(\ref{erbim}), (\ref{lbimr}) and (\ref{rbiml}). 
In this situation the standard notation for the tensor product of
modules, e.g. $A\stackrel{\ot}{_{_R}} A$, would be ambigous. 
For this reason we put marks on both modules as in  
$\asr\ot\atr$, for example, to indicate the
module structures taking part in the tensor products. 

For coproduts of left bialgebroids we use the Sweedler notation in
the form $\gamma_L(a)=a\di \ot a\dii$ and of right bialgebroids
$\gamma_R(a)=a\ui\ot a\uii$.

The Hopf algebroid introduced in \cite{HGD} has both left
and right bialgebroid structures and an antipode relating
the two
\footnote{Actually the original definition in \cite{HGD} was
formulated in terms of the pair $({\cal A}_L,S)$ consisting of a left
bialgebroid and an antipode. It was proven that there exists a right
bialgebroid ${\cal A}_R$ together with which the triple ${\cal
A}=({\cal A}_L,{\cal A}_R,S)$ satisfies the Definition \ref{hgd}
below. The right bialgebroid ${\cal A}_R$ is unique up to an irrelevant
isomorphism.}: 
\bd \lb{hgd} The triple ${\cal A}=({\cal A}_L,{\cal A}_R,S)$ is a {\em
Hopf algebroid} if 
${\cal A}_L=(A,L,s_L,t_L,\gamma_L,\pi_L)$ is a left bialgebroid and
${\cal A}_R=(A,R,s_R,t_R,\gamma_R,\pi_R)$ is a right bialgebroid such
that the base 
rings are related to each other via $R\simeq L^{op}$  and
\be i)\  s_L(L)=t_R(R) \qquad t_L(L)=s_R(R) \lb{schi}\ee
as subrings of $A$,
\bea ii) &(\gamma_L \ot \id_A)\ci \gamma_R &= (\id_A \ot \gamma_R)\ci
\gamma_L \nn
         &(\gamma_R \ot \id_A)\ci \gamma_L &= (\id_A \ot \gamma_L)\ci
\gamma_R
\lb{schii}
\eea
as maps $A\to \atl \ot  {\asl^R} \ot  {\atr} $ and
$A\to \asr \ot  {\atr_L}\ot \asl$,
respectively. 

The map  $S:A\to A$ is  a bijection of additive groups such that 
\be iii)\   S(t_L(l)a t_L(l^{\prime}))=s_L(l^{\prime})S(a) s_L(l) \qquad 
            S(t_R(r^{\prime})a t_R(r))=s_R(r)S(a) s_R(r^{\prime}) 
\lb{defiii} \ee 
for all $l,l^{\prime}\in L$, $r,r^{\prime}\in R$ and $a\in A$. 
The requirement
(\ref{defiii}) makes the expressions $S(a\di)a\dii$ and $a\ui S(a\uii)$
meaningful. The axioms
\be iv)\ S(a\di)a\dii =s_R\ci \pi_R (a) \qquad 
         a\ui S(a\uii)=s_L\ci \pi_L (a) \lb{defiv}
\ee
are required for all $a$ in $A$.
\ed
We emphasize that this notion of Hopf algebroid is different from the
one introduced in \cite{Lu}.

\bd A {\em Hopf algebroid homomorphism } $({\cal A}_L,{\cal A}_R,S)\to 
({\cal A}^{\prime}_L,{\cal A}^{\prime}_R,S^{\prime})$ is a left
bialgebroid homomorphism  $(\Phi,\phi):{\cal A}_L\to {\cal
A}_L^{\prime}$. 

Since it is proven in \cite{HGD} Proposition 4.3 that  both
$(S,\pi_R\ci s_L)$ and $(S^{-1},\pi_R\ci t_L)$ are left bialgebroid
isomorphisms ${\cal A}_L\to ({\cal A}_R)^{op}_{cop}$, for a Hopf
algebroid homomorphism $(\Phi,\phi) : ({\cal A}_L,{\cal A}_R,S)\to 
({\cal A}^{\prime}_L,{\cal A}^{\prime}_R,S^{\prime})$ both
$(S^{\prime}\ci \Phi\ci 
S^{-1}, \pi_R^{\prime}\ci s_L^{\prime}\ci \phi\ci \pi_L\ci t_R)$ and
$( S^{\prime -1}\ci \Phi\ci S,\pi_R^{\prime}\ci t_L^{\prime}\ci
\phi\ci \pi_L\ci s_R)$ are right bialgebroid homomorphisms $:{\cal
A}_R\to {\cal A}_R^{\prime}$.

A {\em Hopf algebroid isomorphism} is a Hopf algebroid homomorphism
$(\Phi,\phi)$ such that both maps $\Phi$ and $\phi$ are bijective.

A Hopf algebroid homomorphism $(\Phi,\phi): ({\cal A}_L,{\cal A}_R,S)\to 
({\cal A}^{\prime}_L,{\cal A}^{\prime}_R,S^{\prime})$ is {\em strict}
if $\Phi\ci S=S^{\prime}\ci \Phi$. 
\ed
The existence of non-strict isomorphisms of Hopf algebroids is a new
feature compared to (weak) Hopf algebras.
\bigskip


In the rest of this subsection let ${\cal A}=({\cal A}_L,{\cal
A}_R,S)$ be a Hopf algebroid where ${\cal
A}_L=(A,L,s_L,t_L,\gamma_L,\pi_L)$ is the left bialgebroid
and ${\cal A}_R=(A,R,s_R,t_R,\gamma_R,\pi_R)$ is the right bialgebroid
underlying ${\cal A}$. 
The left and right integrals in a Hopf algebroid are introduced as the
invariants of the left and right regular modules, respectively:

\bd 
The {\em left integrals} are the elements of the right ideal:
\[ {\cal I}^L({\cal A})\colon = \{ \ell\in A \vert a\ell=s_L\circ
\pi_L(a)\ell \quad \forall a\in A \} .  \]
The {\em right integrals} are the elements of the left ideal:
\[ {\cal I}^R({\cal A})\colon = \{ \Upsilon\in A \vert \Upsilon a=\Upsilon s_R\circ
\pi_R(a)\quad  \forall a\in A\}. \]
\ed
The following lemma is cited from \cite{HGD} Lemma 5.2 and it will
be of importance in the considerations of this paper:
\bl \lb{lemint}
The following characterizations of right/left integrals are equivalent:
\[ \begin{array}{rlcl}
i)&\re\in  {\cal I}^R({\cal A})&\qquad / \qquad&\ell\in \ila\nn
ii)&\re a= \re t_R\circ \pi_R(a) &\qquad / \qquad& a\ell=t_L\ci
\pi_L(a)\ell \quad {for\ all\ } a\in A\nn
iii)& S(\re)\in {\cal I}^L({\cal A})&\qquad / \qquad&S(\ell)\in{\cal I}^R(A)\nn
iv)& S^{-1} (\re)\in {\cal I}^L({\cal A})&\qquad /
\qquad&S^{-1}(\ell)\in{\cal I}^R(A) \nn
v) & \re_{(1)} a \otimes \re_{(2)} = \re_{(1)} \otimes \re_{(2)}S(a) 
)&\qquad / \qquad&\ell\ui\  \ot\  a\ell \uii=S(a)\ell\ui\ \ot\
\ell\uii \nn
&{as\  elements\  of\ } A_L \otimes {_LA} &\qquad \ \qquad&
{as\  elements\  of\ } A^R \otimes {^R A} 
\qquad {\ for\ all\ } a\in A.
\nonumber\end{array}\]
\el
The  left integral $\ell$ was called {\em non-degenerate} in
\cite{HGD} if the maps
\bea \fsr:\asrd\ \to& A \qquad \phisr& \mapsto \ \phisr\ru \ell \quad{\rm
and}\nn
     \ftr\ :\atrd\ \to& A \qquad \phitr& \ \mapsto \ \phitr\rd \ell 
\nonumber\eea
are bijective. It follows that introducing $\lsr\colon
=\fsr^{-1}(1_A)$ 
the inverses of the maps ${\ell_R}$ and ${_R \ell}$ are
written as 
\[ {\ell_R}^{-1}(a) = \lsr\lu S(a) \qquad
   {_R\ell }^{-1}(a)= \lsr\ci S\ld S^{-1}(a). \]
Analogously, the  right integral $\re$ is {\em non-degenerate} if the
maps
\bea {_L\re}:\asld\ \to &A \qquad \phisl &\mapsto\ \re \ld \phisl
\quad {\rm and}\nn
     {\re_L}:\atld\ \to &A \qquad \phitl&\mapsto\ \re\lu \phitl
\nonumber\eea
are bijective. Introducing ${_*\rho}\colon = {_L\re}^{-1}(1_A)$
the inverses of the maps ${_L \re}$ and  ${\re _L}$ are
written as 
\[ {_L \re}^{-1}(a) = S(a)\rd {_*\rho} \qquad
   {\re _L}^{-1}(a) = S^{-1}(a)\ru ({_* \rho}\ci S). \]
It is proven in \cite{HGD} Theorem 5.5 that the existence of a
non-degenerate left integral $\ell$ in the Hopf algebroid ${\cal A}$
implies that all ring extensions 
$s_L:L\to A$, $t_L:L^{op}\to A$, $s_R:R\to A$ and $t_R:R^{op}\to A$
are  Frobenius extensions. In particular the modules $\asl$, $\atl$,
$\asr$ and $\atr$ are finitely generated projective hence all rings
$\asld$, 
$\atld$ $\asrd$ and $\atrd$ carry bialgebroid structures. The
right bialgebroids $\asld_R$ and ${\cal A}_{* R}$ and
also the left bialgebroids 
$\asrd_L$ and $\atrd_L$ were shown to be isomorphic. Furthermore, the left
bialgebroid $\asrd_L$ and
the right bialgebroid $\asld_R$  turned out to be
anti-isomorphic. 

What is more, fixing a non-degenerate left integral $\ell$ in ${\cal
A}$, one can construct an ($\ell$-dependent) Hopf algebroid structure
$\asrd_{\ell}$ on the ring $\asrd$ with two sided non-degenerate
integral ${\ell_R^{-1}}(1_A)$.

The Hopf algebroids $\asrd_{\ell}$ and $\asrd_{\ell^{\prime}}$ --
corresponding to different choices of the non-degenerate left integral
-- are isomorphic but not strictly isomorphic.

Assuming the existence of a {\em two sided} (i.e. both left and right)
non-degenerate integral $i$
\footnote{For a two sided integral $i$ the non-degeneracy as a left
and as a right integral are equivalent.} 
, the 
Hopf algebroid $\asrd_i$ can be interpreted as the dual
Hopf algebroid of ${\cal A}$ in the following sense:
\bt \cite{HGD} \lb{dualthm} 
Let ${\cal A}$ be a Hopf algebroid with two sided
non-degenerate integral $i$. Then the following data define a Hopf
algebroid $\asrd_i$: 
The left bialgebroid over the base $R$
and the right bialgebroid over the base $L$ on the ring $\asrd$ given by 
\[\begin{array}{ll}
     s^*_L(r)(a)=r\pi_R(a)&s^{*}_R(l)(a) = \pi_R(s_L(l) a)\nn
     t^*_L(r)(a)=\pi_R(s_R(r)a)&t^*_R(l) (a) = \lsr (S(i t_L(l)) a) \nn
     \gamma^*_L(\phisr)=\phisr\lu i\ui \ot {i_R}^{-1}(i\uii)\qquad
     &\gamma^*_R(\phisr) =\phisr \lu S^2(i\uii) \ot{i_R}^{-1}(i\ui)\nn
     \pi^*_L(\phisr)=\phisr(1_A)
     &\pi^*_R(\phisr)=\pi_L\ci s_R\ci (\lsr \phisr) (i)
\end{array}\]
where $\lsr={i_R}^{-1}(1_A)$, and the antipode $S^* \equiv {i_R}^{-1}\ci S \ci
{i_R} :\asrd\to \asrd$. By iterating the construction the
second dual  Hopf algebroid $(\asrd_i)^*_{i_R^{-1}(1_A)}$ is strictly
isomorphic to ${\cal A}$. 
\et
The above Theorem \ref{dualthm} is a variant of \cite{HGD} Theorem
5.17.

\section{Hopf algebroid symmetry of abstract D2 Frobenius extensions}
\setc{0}
\lb{bic}

\subsection{Some harmonic analysis}
\lb{harm}


We adopt a graphical notation for
bicategories using planar diagrams similar to that has been used e.g. by Yetter
for sovereign monoidal categories in \cite{Yetter}. 
The only difference is that now also the planar regions carry labels, the
objects of the bicategory. As an experiment, we will also employ a notation
for the coherence isomorphisms at the price of introducing some metrical
information into the diagrams which, of course, destroys the topological
nature of 2-categorical diagrams. 

Let $(\c,s_0,t_0,s_1,t_1,\circ,\times,\l,\r,\a)$ be a
bicategory.  We use the following graphical rules.
For a 1-cell $a$ we draw a vertical line. The area right to
the line is labelled by the source 0-cell $s_0(a)$ and the area
left to the line by the target $t_0(a)$. For the 2-cell $x$ we
draw a box with upper `leg' its source $s_1(x)$ and lower `leg' the
target $t_1(x)$. The vertical composition $x\ci y$
is represented by a picture in which $x$ is placed under $y$.
The horizontal composition $x\x y$ is represented by a picture in which
$x$ is placed left to $y$. Since the horizontal composition is not
strictly associative we take care about the horizontal distances
between two lines or boxes. The successive bracketing is represented
by growing distances. For the coherence natural isomorphisms we
use the graphical notation

\begin{figure}[h]
\psfrag{alpha}{\Large$\a_{a,b,c}=$}
\psfrag{l}{\Large$\l_a=$}
\psfrag{r}{\Large$\r_a=$}
\psfrag{t}{\Large$t_0(a)$}
\psfrag{s}{\Large$s_0(a)$}
\psfrag{a}{\Large$a$}
\psfrag{b}{\Large$b$}
\psfrag{c}{\Large$c$}
\begin{center}
{\resizebox*{8cm}{!}{\includegraphics{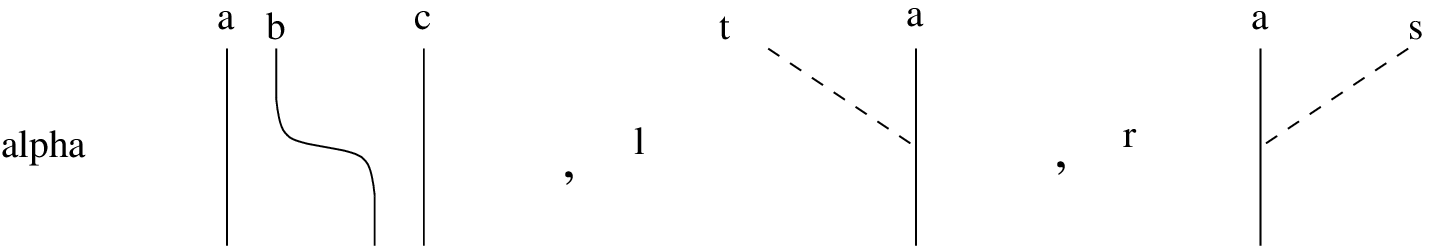}}
}
\end{center}
\end{figure}
\noindent
where the dashed line emphasizes that we have a 0-cell here. 

Let $\iota$ be a
Frobenius 1-morphism of $\c$
-- that is suppose that $\iota$ has a two sided dual ${\bar{\iota}}$.
The Frobenius property  means the existence of 2-morphisms 
\bea &\ev_R\in \c^2(\iota\times {\bar{\iota}},t_0(\i)) \qquad
     &\coev_R\in \c^2(s_0(\i),{\bar{\iota}}\times \iota)\nn
     &\ev_L\in \c^2({\bar{\iota}}\times {{\iota}},s_0(\i)) \qquad
     &\coev_R\in \c^2(t_0(\i),{{\iota}}\times {\bar{\iota}})\lb{rigint}\eea
satisfying the relations
\bea \r_{\iota}\ci (\iota\times \ev_L)\ci \a_{\i,\ib,\i}\ci
(\coev_L\times \iota)\ci  \l_{\iota}^{-1}&=&\iota \nn
\l_{\bar{\iota}}\ci (\ev_L\times {\bar{\iota}})\ci
\a_{\ib,\i,\ib}^{-1} \ci (\ib \x \coev_L)\ci \r_{\ib}^{-1}&=&\ib\nn
\l_{\i}\ci (\ev_R\x \i)\ci \a_{\i,\ib,\i}^{-1}\ci (\i\x \coev_R)\ci
\r_{\i}^{-1} &=&\i\nn
\r_{\ib}\ci (\ib \x \ev_R)\ci \a_{\ib,\i,\ib}\ci(\coev_R\x \ib)\ci
\l_{\ib}^{-1}&=&\ib. \lb{rigrel}\eea
In the graphical notation we draw a vertical line directed
downwards for $\i$, upwards for $\ib$ and 
\vspace{0.3cm}
\begin{center}
{\resizebox*{10cm}{!}{\includegraphics{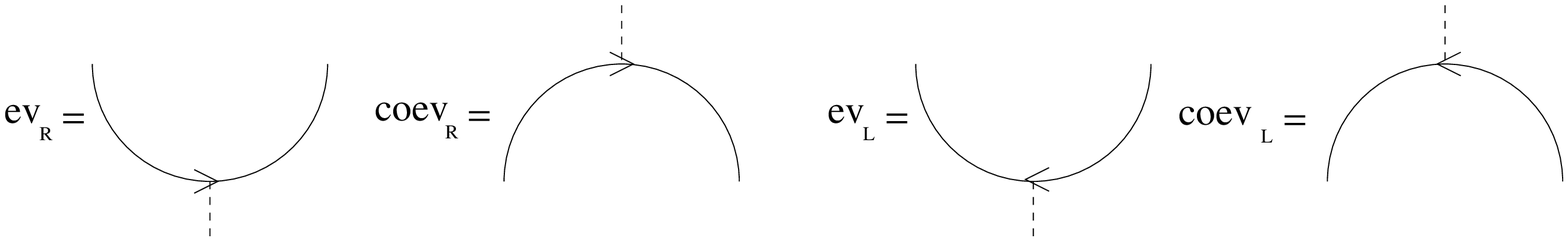}} 
\vspace{0.3cm}}
\end{center}
In this language the relations (\ref{rigrel}) read as
\vspace{0.3cm}
\begin{center}
{\resizebox*{2cm}{!}{\includegraphics{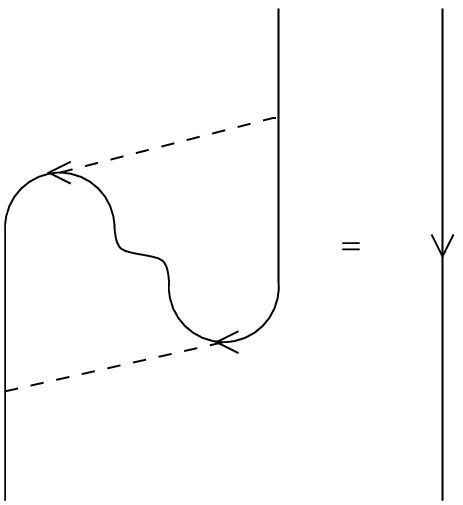}} ~~~~~~~~~~~~~~~~~~
\resizebox*{2cm}{!}{\includegraphics{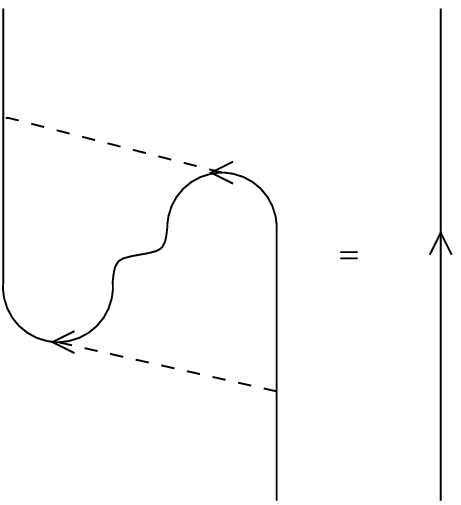}}   ~~~~~~~~~~~~~~~~~~
\resizebox*{2cm}{!}{\includegraphics{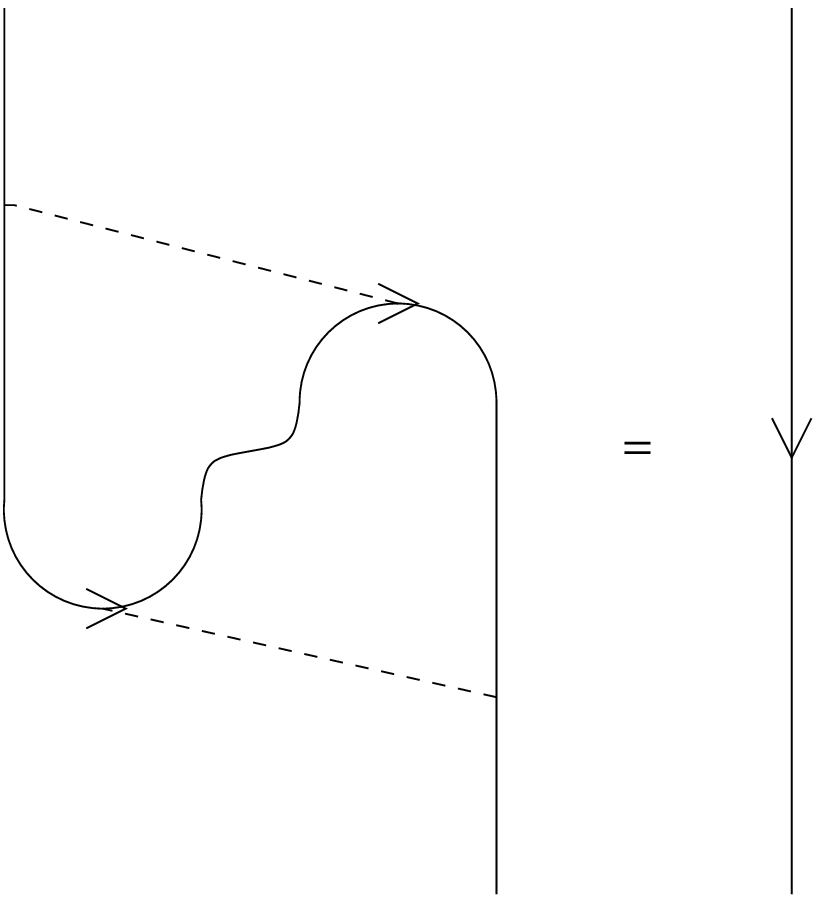}}  ~~~~~~~~~~~~~~~~~~
\resizebox*{2cm}{!}{\includegraphics{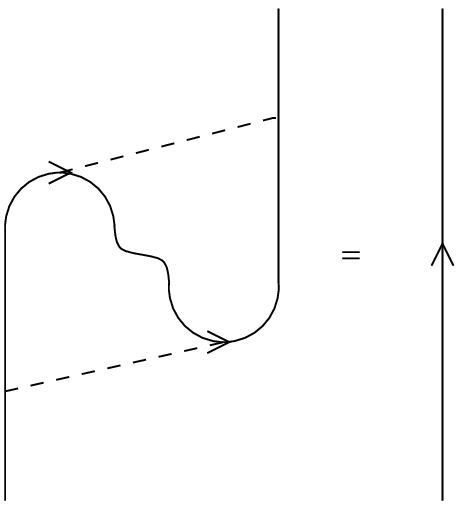}} \vspace{0.3cm}}
\end{center}

It follows from coherence that the meanderings of
the lines as well as the dashed lines may be forgotten if equality of two
2-cells is to be proven graphically. 

\bigskip


From now on we
assume that the bicategory $\c$ is enriched over ${\bf Ab}$, i.e., the 
hom sets are Abelian groups, written additively, and both the vertical and the
horizontal  compositions are group homomorphisms. In this situation
both isomorphic additive groups  ${\tt A}\colon =
\c^2(\i\x\ib,\i\x\ib)$ and  ${\tt B}\colon = \c^2(\ib\x\i,\ib\x\i)$
carry two ring structures. The multiplication in $A\colon =({\tt
A},\ci)$ and $B\colon =({\tt B},\ci)$ is given by the vertical composition
$\ci$ of $\c$ while ${\hat A}\colon =({\tt A},\ast)$ and ${\hat
B}\colon = ({\tt B},\ast)$ can be defined -- by fixing the 2-morphisms
(\ref{rigint}) -- as  
\bea a_1\ast a_2 &=& (\r_{\i}\x \ib)\ci 
[(\i\x {\ev_L})\x \ib]\ci
(\a_{\i,\ib,\i} \x\ib)\ci 
\a^{-1}_{\i\x\ib,\i,\ib} \ci 
(a_1\x a_2)\ci \nn
&\ci&\a_{\i\x\ib,\i,\ib}  \ci 
(\a^{-1}_{\i,\ib,\i} \x\ib)\ci 
[(\i\x {\coev_R})\x \ib]\ci
 (\r^{-1}_{\i}\x \ib)\lb{conva}\\
b_1\ast b_2 &=& (\ib\x \l_{\i})\ci
[\ib\x({\ev_R}\x\i)]\ci
(\ib\x\a^{-1}_{\i,\ib,\i})\ci
\a_{\ib,\i,\ib\x\i}\ci
(b_1\x b_2)\ci \nn
&\ci&\a^{-1}_{\ib,\i,\ib\x\i}\ci
(\ib\x\a_{\i,\ib,\i})\ci
[\ib\x({\coev_L}\x\i)]\ci
 (\ib\x \l^{-1}_{\i}).\lb{convb}
\eea
In the graphical language

\begin{figure}[h]
\psfrag{lhs}{\Huge$a_1\ast a_2 =$}
\psfrag{a1}{\Huge$a_1$}
\psfrag{a2}{\Huge$a_2$}
\psfrag{lhsb}{\Huge$b_1\ast b_2 =$}
\psfrag{b1}{\Huge$b_1$}
\psfrag{b2}{\Huge$b_2$}
\begin{center}
{\resizebox*{3cm}{!}{\includegraphics{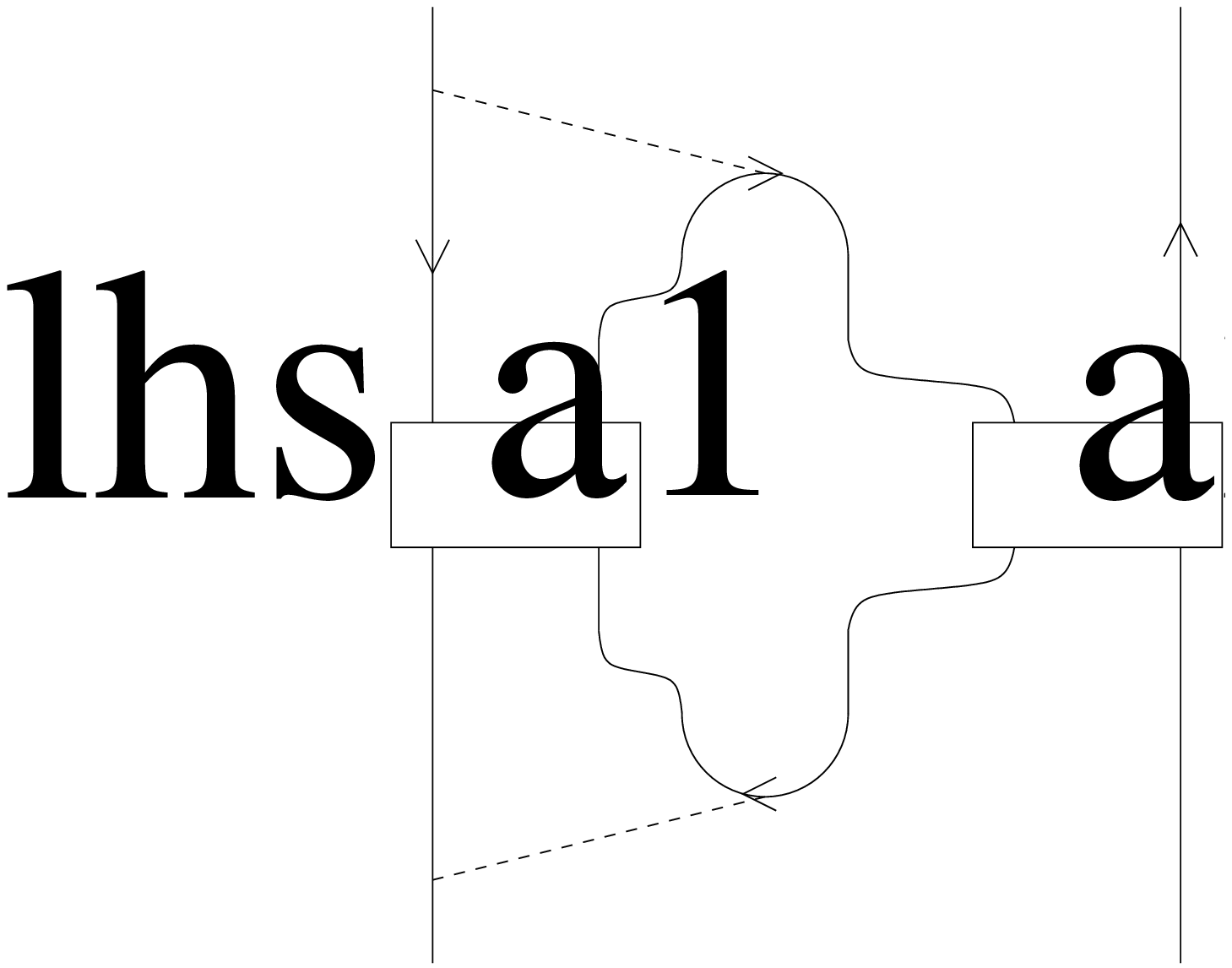}}
~~~~~~~~~~~~~~~~~~~~~~~~~~
\resizebox*{3cm}{!}{\includegraphics{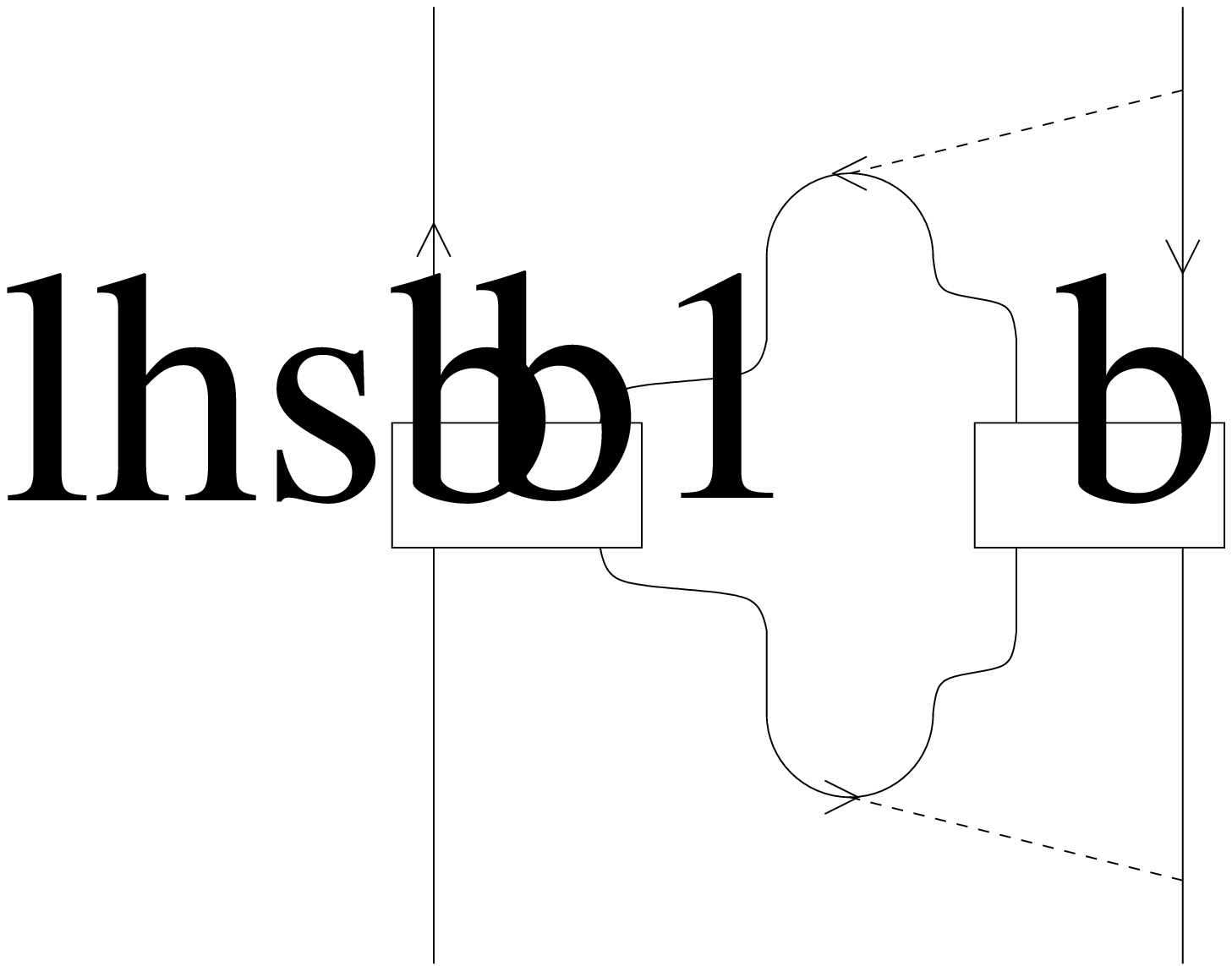}}}
\end{center}
\end{figure}

By coherence $\hat A$ and $\hat B$ are associative with units
$i_A\colon ={\coev_L}\ci {\ev_R}$ and ${i_B}\colon = {\coev_R}\ci
{\ev_L}$, respectively. 
This way we have four associative unital rings 
$A$, $B$, ${\hat A}$ and ${\hat B}$. We claim
that $A$ is isomorphic to ${\hat B}$ and $B$ is isomorphic to $\hat
A$. Let us define the Fourier transformations 
\bea \f: {\tt A}&\to\  {\tt B} \qquad 
a\ \mapsto & \r_{\ib\x\i}\ci
[(\ib\x\i)\x {\ev_L}]\ci
\a_{\ib\x\i,\ib,\i}\ci
(\a^{-1}_{\ib,\i,\ib}\x \i)\ci
[(\ib\x a)\x \i]\ci\nn
&&(\a_{\ib,\i,\ib}\x \i)\ci
\a^{-1}_{\ib\x\i,\ib,\i}\ci
[{\coev_R}\x (\ib\x\i)]\ci
\l^{-1}_{\ib\x\i}\nn
\fdot: {\tt A}&\to\  {\tt B} \qquad 
a\ \mapsto & \l_{\ib\x\i}\ci
[{\ev_L}\x (\ib\x\i)]\ci
\a^{-1}_{\ib,\i,\ib\x\i}\ci
(\ib\x \a_{\i,\ib,\i})\ci
[\ib\x(a\x\i)]\ci\nn
&&
(\ib\x \a^{-1}_{\i,\ib,\i})\ci
\a_{\ib,\i,\ib\x\i}\ci
[(\ib\x \i)\x {\coev_R}]\ci
\r^{-1}_{\ib\x\i}
\eea
In picture:

\begin{figure}[h]
\psfrag{a}{\Huge$a$}
\begin{center}
{\resizebox*{4cm}{!}{\includegraphics{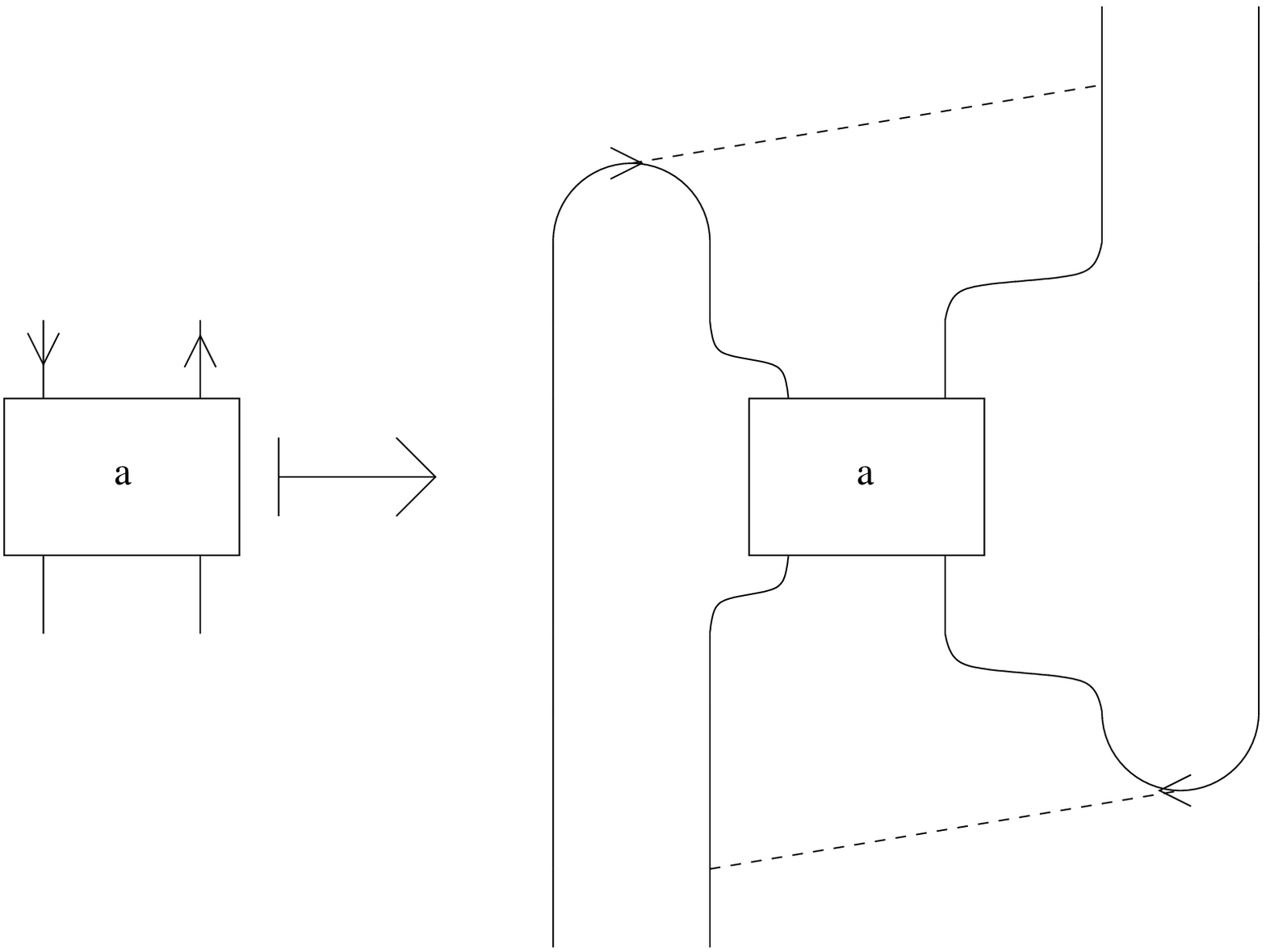}}
~~~~~~~~~~~~~~~~~~~~~~~~~~
\resizebox*{4cm}{!}{\includegraphics{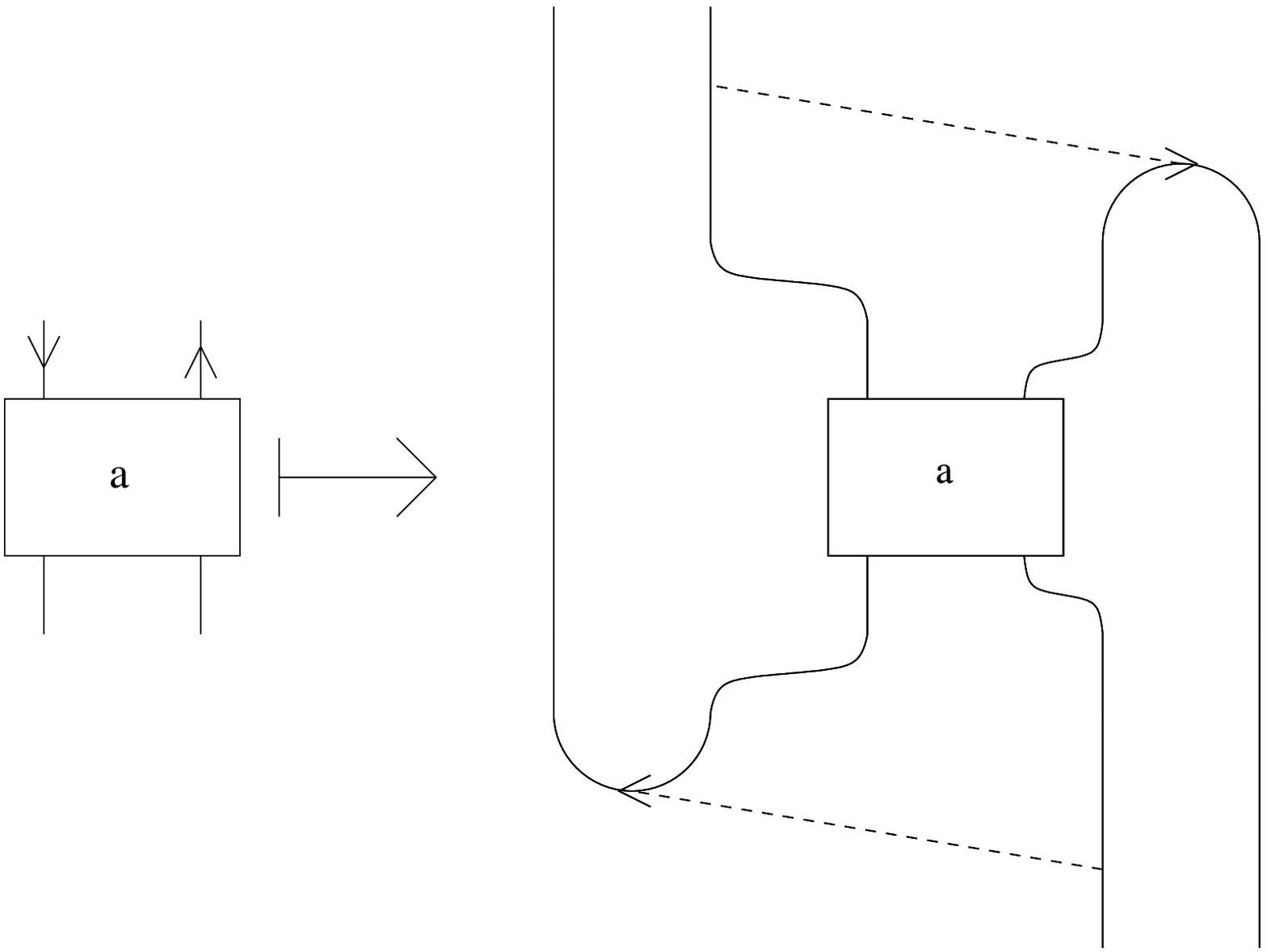}}}
\end{center}
\end{figure}

\noindent
Both $\f$ and $\fdot$ are bijections with inverses
\bea \f^{-1}: {\tt B}&\to \ {\tt A} \qquad
b\ \mapsto & \l_{\i\x\ib}\ci
[{\ev_R}\x(\i\x\ib)]\ci
\a^{-1}_{\i,\ib,\i\x\ib}\ci
(\i\x\a_{\ib,\i,\ib})\ci
[\i\x(b\x\ib)]\ci\nn
&&(\i\x\a^{-1}_{\ib,\i,\ib})\ci
\a_{\i,\ib,\i\x\ib}\ci
[(\i\x\ib)]\x{\coev_L}]\ci
\r^{-1}_{\i\x\ib}\nn
\fdot^{-1}: {\tt B}&\to \ {\tt A} \qquad
b\ \mapsto & \r_{\i\x\ib}\ci
[(\i\x\ib)\x{\ev_R}]\ci
\a_{\i\x\ib,\i,\ib}\ci
(\a^{-1}_{\i,\ib,\i}\x \ib)\ci
[(\i\x b)\x\ib]\ci \nn
&&(\a_{\i,\ib,\i}\x \ib)\ci
\a^{-1}_{\i\x\ib,\i,\ib}\ci
[{\coev_L}\x(\i\x\ib)]\ci
\l^{-1}_{\i\x\ib}
\eea

\noindent
They relate the products $\ci$ and $\ast$ as follows:
\bea \f(a_1\ci a_2)&=& \f(a_2)\ast \f(a_1)\nn
\fdot(a_1\ci a_2)&=& \fdot(a_1)\ast \fdot(a_2)\nn
\f(a_1\ast a_2)&=& \f(a_1)\ci \f(a_2)\nn
\fdot(a_1\ast a_2)&=& \fdot(a_2)\ci \fdot(a_1).\lb{fri}
\eea
The equations (\ref{fri}) imply that the differences of the Fourier
transformations give ring anti-automorphisms 
\be S_A\colon = \fdot^{-1}\ci \f:A^{op}\to A
\quad {\rm and} \quad  
S_B\colon = \f\ci \fdot^{-1}: B^{op}\to B.
\lb{s}\ee 
They are to be the antipodes in
the case when $\iota$ satisfies the D2 condition discussed in the next
subsection.

Let us investigate some bimodule map properties of the maps $S_A$ and
$S_B$. 
Both additive groups ${\tt A}$ and ${\tt B}$ carry four commuting
actions of the anti-isomorphic rings $L\colon = \c^2(\i,\i)$ and
$R\colon = \c^2(\ib,\ib)$ -- given by `composition on the four
legs'. That is we have the ring isomorphisms 
\bea \mu: L&\to\ R^{op}\qquad l&\mapsto\ 
\l_{\ib}\ci ({\ev_L}\x \ib)\ci [(\ib \x l)\x\ib]\ci
\a^{-1}_{\ib,\i,\ib} \ci (\ib\x{\coev_L})\ci \r^{-1}_{\ib} \nn
\nu:  L&\to\ R^{op}\qquad l&\mapsto\ 
\r_{\ib}\ci (\ib\x {\ev_R})\ci [\ib\x(l\x\ib)]\ci\a_{\ib,\i,\ib}\ci
({\coev_R}\x \ib)\ci \l_{\ib}^{-1}
\lb{munu}\eea
which read in the graphical language as

\begin{figure}[h]
\psfrag{mu}{\Huge$\mu(l)=$}
\psfrag{nu}{\Huge$\nu(l)=$}
\psfrag{l}{\Large$l$}
\begin{center}
{\resizebox*{7cm}{!}{\includegraphics{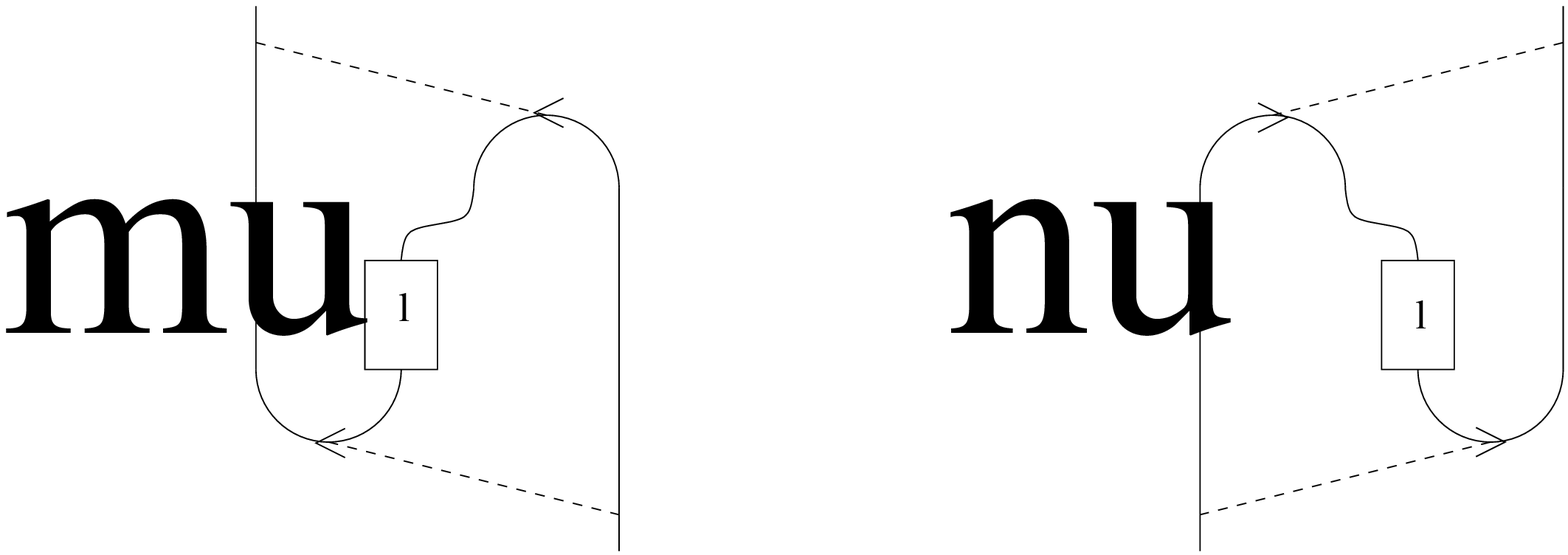}}}
\end{center}
\end{figure}

\noindent
With the help of the maps $\mu$ and $\nu$ we can define the bimodules
${_LA_L}$ and ${^RA^R}$  as follows:
\bea l\cdot a &= (l\x\ib)\ci a\qquad
     a\cdot l &= (\i\x \mu(l))\ci a\nn
     a\cdot r &= a\ci (\i\x r)\qquad
     r\cdot a &= a\ci (\nu^{-1}(r)\x \ib).
\lb{bim}\eea
Since
\bea \lb{convmod}
\f(l\cdot a) \ =& (\ib\x l)\ci \f(a) \qquad
\f(a\cdot l) &=\ \f(a)\ci (\ib\x l)\nn
\f(a\cdot r) \ =& \f(a)\ci (r\x \i) \qquad
\f(r\cdot a) &=\ (r\x\i)\ci \f(a)
\eea
and 
\bea \lb{convmodd}
\fdot(l\cdot a)\ =& \fdot(a)\ci (\mu(l)\x\i)\qquad\quad
\fdot(a\cdot l)&=\  (\mu(l)\x\i)\ci \fdot(a)\nn
\fdot(a\cdot r)\ =& (\ib\x \nu^{-1}(r))\ci \fdot(a)\qquad
\fdot(r\cdot a)&=\ \fdot(a)\ci (\ib\x \nu^{-1}(r))
\eea
the map $S_A$  is compatible with the
bimodule structures (\ref{bim}) that is it is a {\em twisted bimodule
map} \cite{HGD} ${_LA_L}\to{^RA^R}$ and also ${^RA^R}\to {_LA_L}$:
\be S_A(l\cdot a\cdot l^{\prime})= \nu( l^{\prime})\cdot a\cdot
\nu(l)\qquad
     S_A(r\cdot a \cdot r^{\prime})= \mu^{-1}( r^{\prime})\cdot a
\cdot \mu^{-1}(r).
\ee

\subsection{The depth 2 case}
\lb{depth2}


The notion of {\em depth 2} or shortly D2 property in the context of
bicategories was introduced in \cite{Sz}
\footnote{Recall that the terminology of \cite{Sz} is somewhat
different from the later publications. The condition that was called
D2 property in \cite{Sz} we call {\em left} D2 property, as it is
explaned in the Introduction.}: Let $\c$ be an ${\bf Ab}$-enriched
bicategory closed under direct sums of 1-morphisms and possessing
zero 1-morphisms for any pair of objects. 
Such bicategories will be called {\em additive}. Let
$\i$ be a 1-morphism in $\c$ possessing a left dual $\ib$. Then $\i$
is said to satisfy the left D2 condition if $(\i\x\ib)\x\i$ is a direct
summand in a finite direct sum of $\i$'s. In this case -- under the
additional assumption that $s_0(\i)$ is a direct
summand in $\ib\x\i$, which can be relaxed -- it was proven in \cite{Sz}
that for such a 1-morphism $\iota$ the ring $A=
(\c^2(\i\x\ib,\i\x\ib),\ci)$ has a canonical left bialgebroid
structure over the base $L=\c^2(\i,\i)$. 

It is clear that if  $(\ib\x\i)\x\ib$ is a direct
summand in a finite direct sum of $\ib$'s (right D2 condition) then
the ring $B\colon 
= (\c^2(\ib\x\i,\ib\x\i),\ci)$ has a right bialgebroid structure.

In the sequel we will see that if $\i$ is a D2 Frobenius 1-morphism

then so is $\ib$. In this case the  rings $A=(\c^2(\i\x\ib,\i\x\ib),\ci)$
and $B=(\c^2(\ib\x\i,\ib\x\i),\ci)$ carry left as well as right
bialgebroid structures such that together with
the antipodes $S_A$ and $S_B$ in (\ref{s}) they are dual Hopf
algebroids in the sense of \cite{HGD}.  

Throughout the subsection let
$(\c,s_0,t_0,s_1,t_1,\ci,\times,\l,\r,\a)$ be an additive bicategory.
Let $\i$ be a
Frobenius 1-morphism in $\c$. In this case one can reformulate the D2
condition as follows:

\bp The $\i$ satisfies the left D2 condition if and only if there
exists 
an element $\sum_i y_i\ot x_i\in {A^L}\ot{_L A}$
such that 
\be (y_i\x\i)\ci \a_{\i,\ib,\i}^{-1}\ci(\i\x{\coev_R})\ci (\i\x
{\ev_L}) \ci \a_{\i,\ib,\i}\ci(x_i \ot \i)=(\i\x\ib)\x\i \lb{d2}\ee
where $A^L$ is the right $L$ module defined as $a\cdot l \colon = a\ci
(l\x\ib)$. 
In the graphical notation (\ref{d2}) reads as

\begin{figure}[h]
\psfrag{x}{\Huge$x_i$}
\psfrag{y}{\Huge$y_i$}
\begin{center}
{\resizebox*{3cm}{!}{\includegraphics{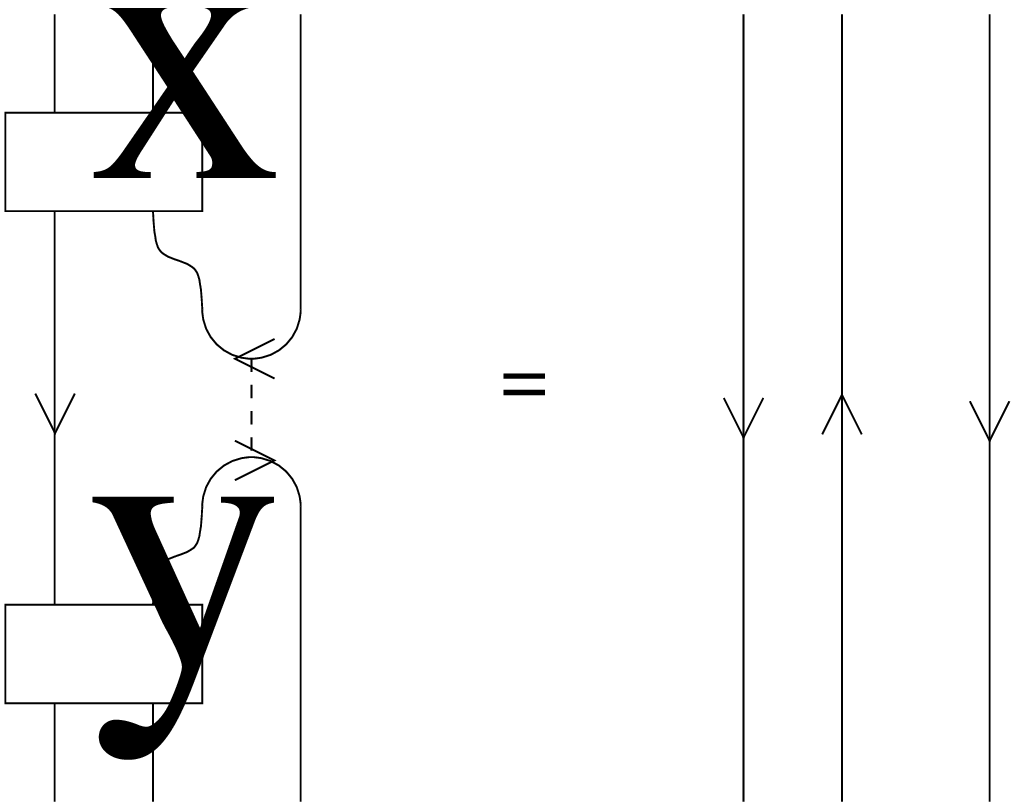}}}
\end{center}
\end{figure}

\noindent
The element $\sum_i y_i\ot x_i\in {A^L}\ot{_L A}$ is
called a {\em D2 quasi-basis} for $\i$.  
\ep
\pr The property that $(\i\x\ib)\x\i$ is a direct summand in a finite
direct sum of $\i$'s means the existence of finite sets of 2-cells
\bea \{ \beta_i\}_{i=1\dots n}&\subset& \c^2((\i\x\ib)\x\i,\i) \nn
     \{ \beta^{\prime}_i\}_{i=1\dots n}&\subset& \c^2(\i,(\i\x\ib)\x\i) 
\eea
satisfying $\sum_i \beta^{\prime}_i\ci  \beta_i=(\i\x\ib)\x\i$. By the
Frobenius property of $\i$ this is further equivalent to the 
existence of the sets
\bea \{x_i\colon = (\beta_i\x\ib)\ci \a^{-1}_{\i\x\ib,\i,\ib}\ci
[(\i\x\ib)\x {\coev_L}]\ci \r^{-1}_{\i\x\ib}\}_{i=1\dots n} &\subset &
\c^2(\i\x\ib,\i\x\ib)\nn 
\{\ y_i\colon = \r_{\i\x\ib}\ci [(\i\x\ib)\x{\ev_R}]\ci
\a_{\i\x\ib,\i,\ib}\ci (\beta^{\prime}_i\x \ib)\ \  \}_{i=1\dots n}
&\subset &  \c^2(\i\x\ib,\i\x\ib)
\eea

\begin{figure}[h]
\psfrag{x}{\Huge$x_i=$}
\psfrag{y}{\Huge$y_i=$}
\psfrag{bi}{\Large$\beta_i$}
\psfrag{bip}{\Large$\beta_i^{\prime}$}
\begin{center}
{\resizebox*{2.5cm}{!}{\includegraphics{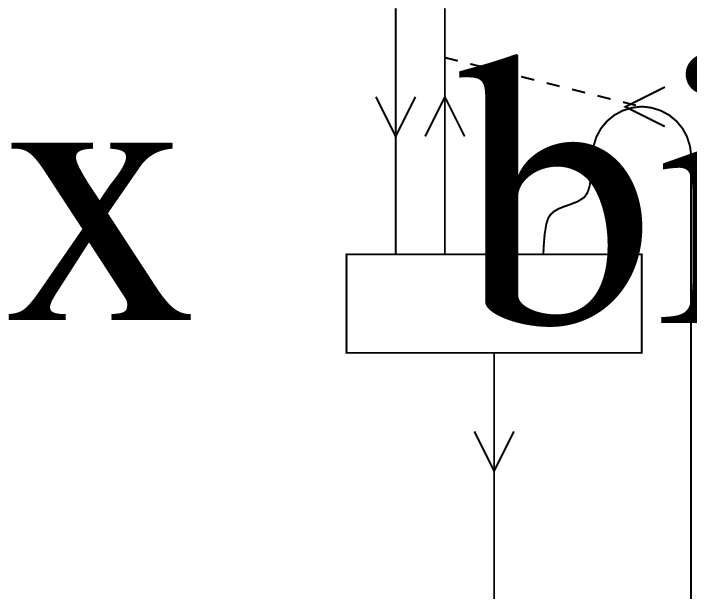}}
~~~~~~~~~~~~~~~~~~~~~~
\resizebox*{2.5cm}{!}{\includegraphics{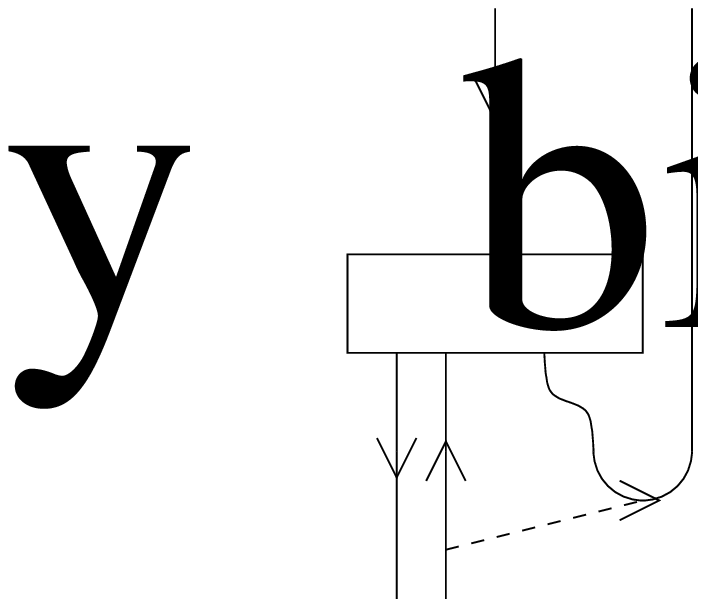}}}
\end{center}
\end{figure}

\noindent
satisfying (\ref{d2}). 
\qed

If $\sum_i y_i\ot x_i$ is a D2 quasi-basis for $\i$ then $\sum_i \f\ci
\fdot^{-1}\ci \f(x_i)\ot \f(y_i) \in B^R \stac{R} {_R B}$ --  where ${_R
B}$ is the left $R$-module defined as $r\cdot  b\colon = (r\x\i)\ci b$
and $B^R$ is the right $R$-module $b\cdot r=b\ci (r \x \i)$
-- is easily checked to be a D2 quasi-basis for $\ib$. Since, owing
to the Frobenius property of $\i$, the 
right D2 condition on $\i$ coincides with the left D2 condition on
$\ib$ we conclude that $\i$ is D2 if and only if it is left D2 and if
and only if $\ib$ is D2.


From now on let $\i$ be a D2 Frobenius 1-morphism in $\c$ with D2
quasi-basis $y_i\ot x_i$. We omit the summation symbol for summing
over the D2 quasi-basis. 
The next Proposition shows that $y_i\ot x_i$ is really a
quasi-basis in the sense of \cite{W}: 
\bp \lb{slfrob} 
The map 
\be \lb{sl}
    \sla:L\to \ A\qquad l\mapsto \ l\x\ib \ee
is a Frobenius extension.    
\ep

\pr We construct the Frobenius map 
\be \lb{phil}\phi_L:A\to L\qquad
a \mapsto \r_{\i}\ci (\i\x {\ev_L})\ci \a_{\i,\ib,\i}\ci (a\x
\i)\ci \a^{-1}_{\i,\ib,\i}\ci (\i\x {\coev_R})\ci \r^{-1}_{\i}\ee

\begin{figure}[h]
\psfrag{lhs}{\Huge$\phi_L(a)=$}
\psfrag{a}{\Large$a$}
\begin{center}
{\resizebox*{3.5cm}{!}{\includegraphics{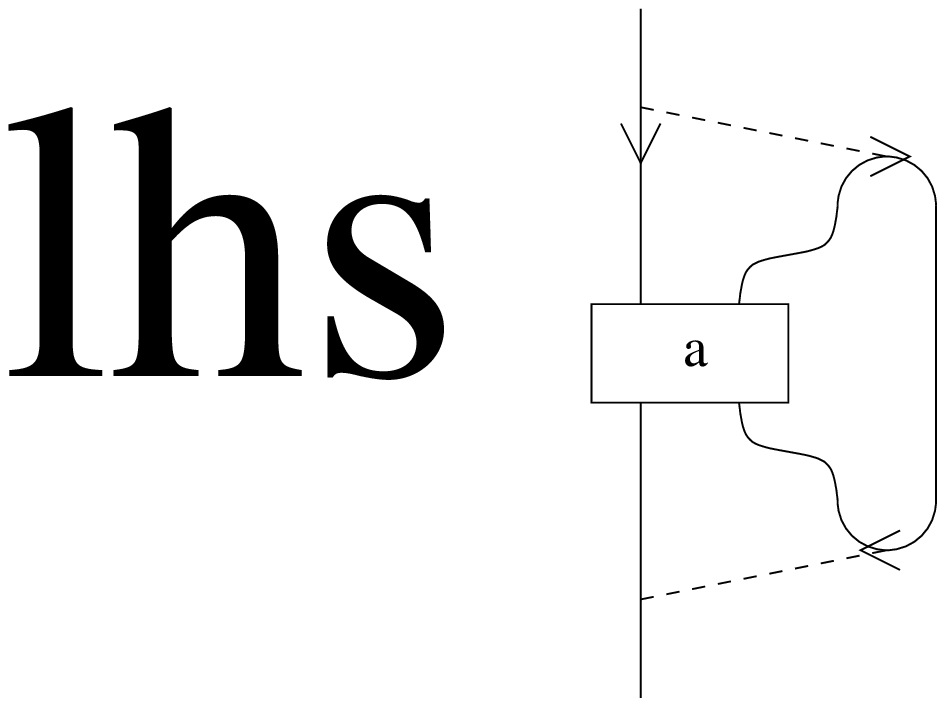}}}
\end{center}
\end{figure}

\noindent
and show that it has the quasi basis $y_i\ot x_i$. Making use of the
coherence axioms in $\c$ 
and the relations (\ref{rigrel}) and (\ref{d2}) we have
\bea \sla\ci \phi_L(a\ci y_i)\ci x_i &=&
(\r_{\i}\x \ib)\ci [(\i \x {\ev_L})\x \ib] \ci (\a_{\i,\ib,\i}\x
\ib) \ci [(a\x\i)\x \ib]\ci  [(y_i \x\i)\x \ib]\ci
(\a^{-1}_{\i,\ib,\i}\x \ib) \ci \nn
&&[(\i\x {\coev_R})\x \ib]\ci
(\r_{\i}^{-1} \x \ib)\ci x_i=\nn
&=&(\r_{\i}\x \ib)\ci [(\i \x {\ev_L})\x \ib] \ci (\a_{\i,\ib,\i}\x
\ib) \ci [(a\x\i)\x \ib]\ci  [(y_i \x\i)\x \ib]\ci
(\a^{-1}_{\i,\ib,\i}\x \ib) \ci \nn
&&[(\i\x {\coev_R})\x \ib]\ci
(\r_{\i}^{-1} \x \ib)\ci (\i\x \l_{\ib})\ci [\i \x ({\ev_L}\x \ib)]\ci 
(\i \x \a_{\ib,\i,\ib}^{-1})\ci \nn
&&[\i\x(\ib \x{\coev_L})]\ci (\i\x \r_{\ib}^{-1}) \ci x_i=\nn
&=&(\r_{\i}\x \ib)\ci [(\i \x {\ev_L})\x \ib] \ci (\a_{\i,\ib,\i}\x
\ib) \ci [(a\x\i)\x \ib]\ci  [(y_i \x\i)\x \ib]\ci
(\a^{-1}_{\i,\ib,\i}\x \ib) \ci \nn 
&&[(\i\x {\coev_R})\x \ib]\ci
[(\i\x {\ev_L})\x\ib]\ci (\a_{\i,\ib,\i}\x \ib)\ci [(x_i\x \i)\x
\ib] \ci \a_{\i\x\ib,\i,\ib}^{-1}\ci \nn
&&[(\i\x\ib)\x {\coev_L}]\ci \r_{\i\x\ib}^{-1} =
a.
\nonumber\eea
Analogously,

\begin{figure}[h]
\psfrag{lhs}{\Huge$y_i\ci \sla\ci \phi_L(x_i\ci a)\ =$}
\psfrag{a}{\Huge$a$}
\psfrag{x}{\Huge$x_i$}
\psfrag{y}{\Huge$y_i$}
\psfrag{=}{\Huge$=$}
\begin{center}
{\resizebox*{6cm}{!}{\includegraphics{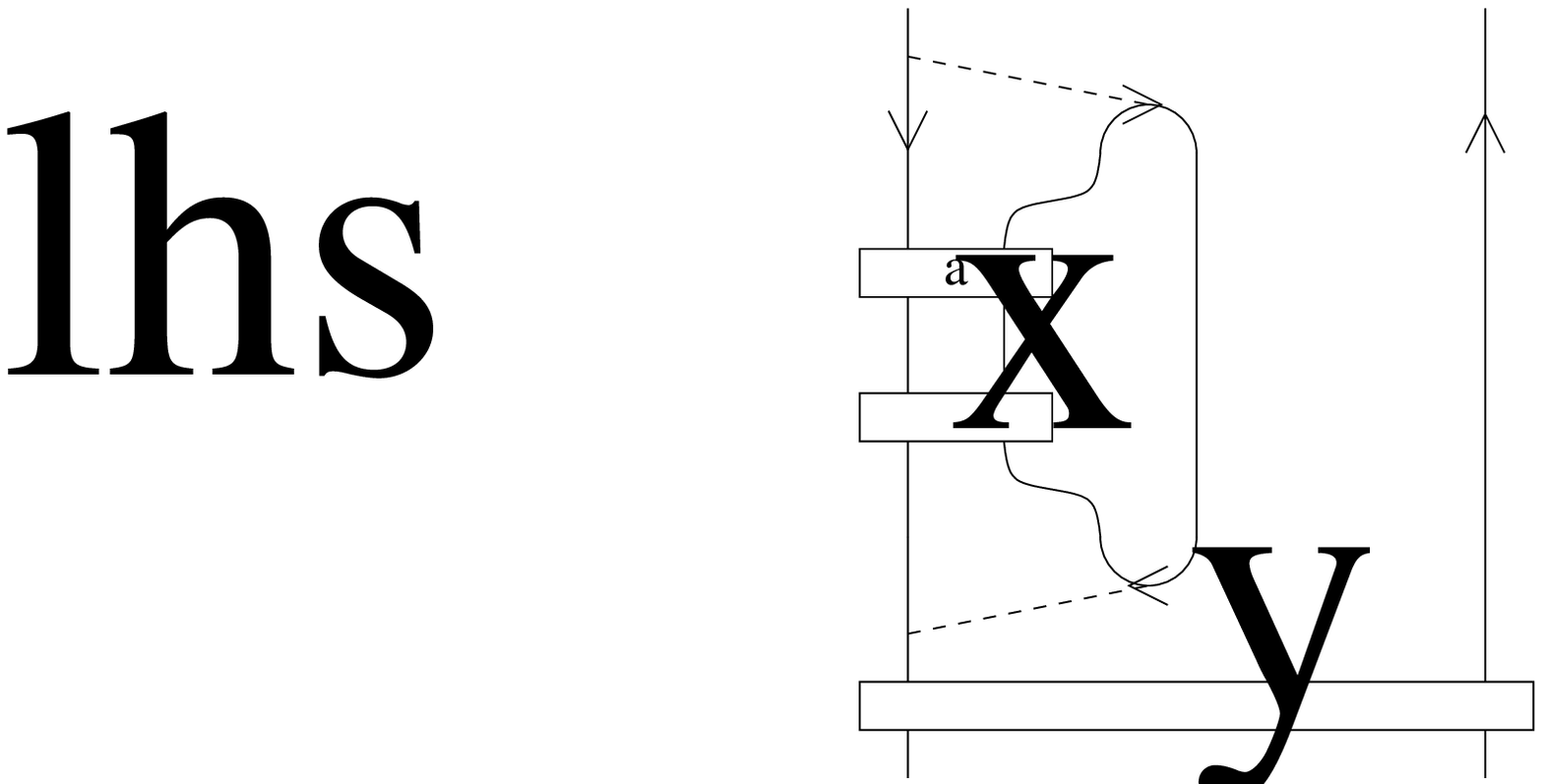}}
~~~~
\resizebox*{2.5cm}{!}{\includegraphics{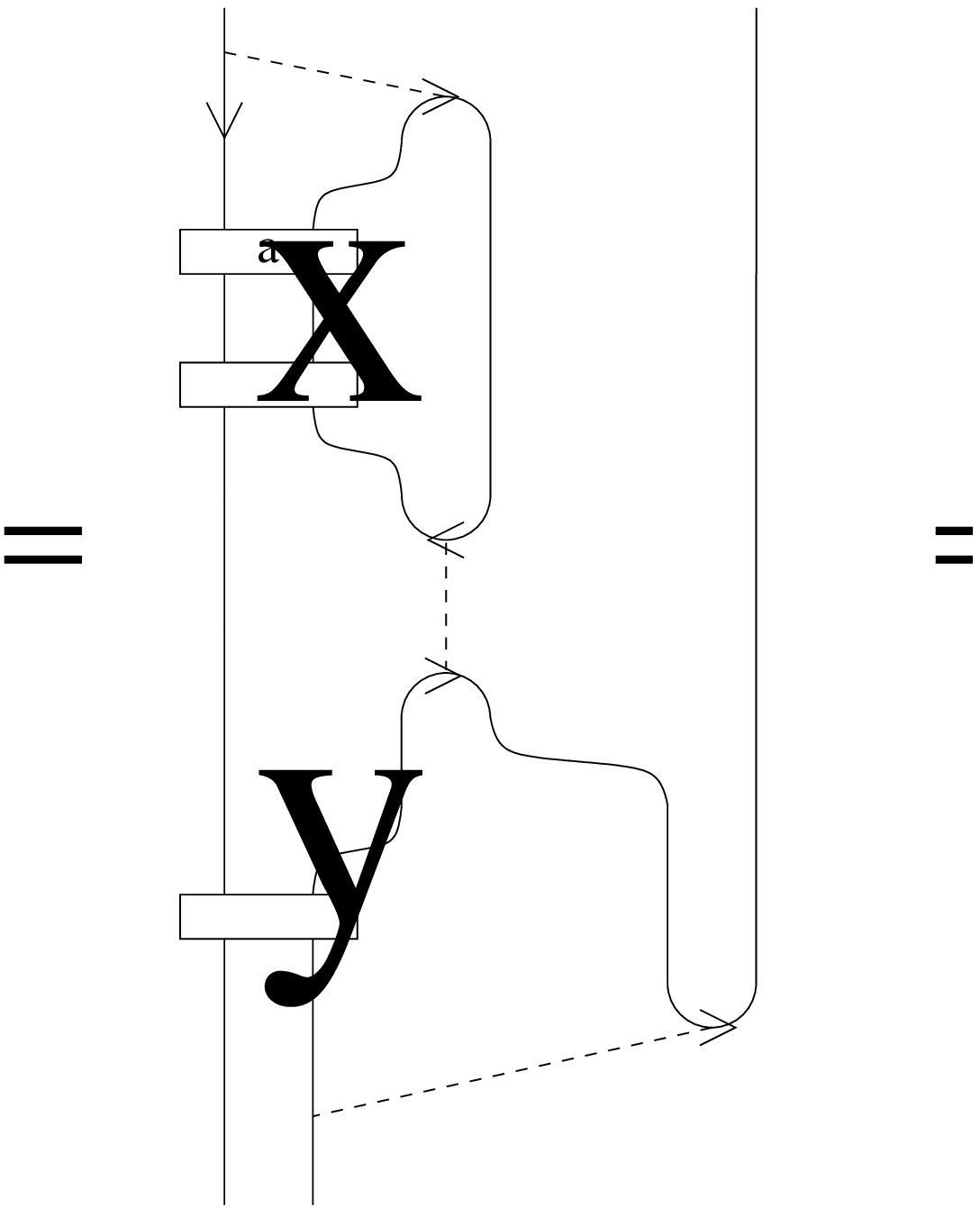}}
~~~~
\resizebox*{4cm}{!}{\includegraphics{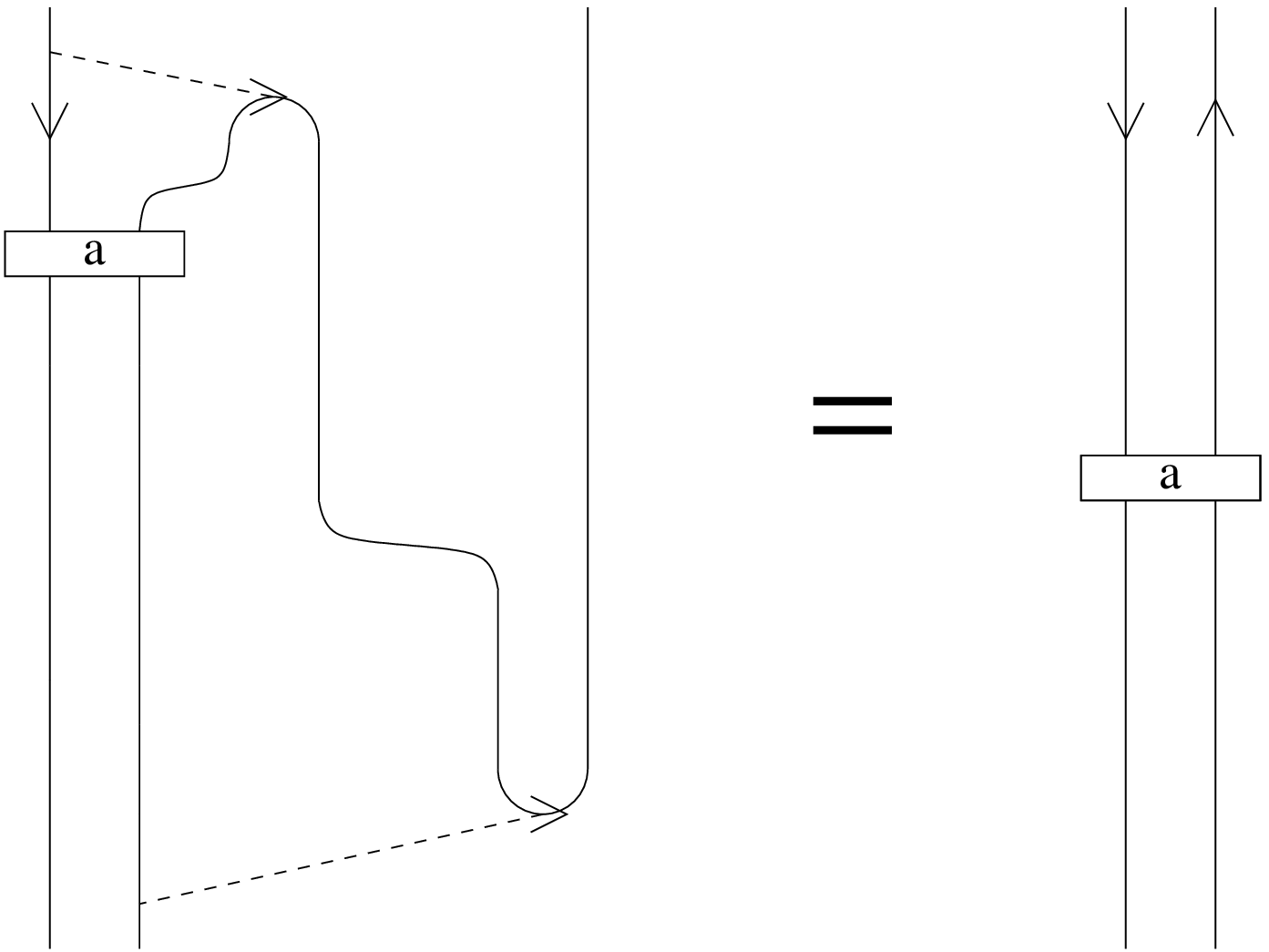}}}
\end{center}
\end{figure}

\qed

\noindent
The following lemma collects some useful properties of the quasi-basis
$y_i\ot x_i$:
\bl For any elements $a,a_1,a_2,a_3$ in $A$ the following hold true:
\bea i)&  a\ci y_i\ot x_i=y_i\ot x_i \ci a \quad {\rm as\  elements\
of\ } {A^L}\ot {_L A} \lb{sep}\\
ii)&  y_i \ot x_i \ast a = y_i \ast S_A(a) \ot x_i \quad {\rm as\  elements\
of\ } {A^L}\ot {_L A} \lb{lem326}\\ 
iii)&  a_1\ast (a_2\ci a_3)\ =\ [(a_1\ci y_i)\ast a_2]\ci (x_i \ast
a_3)\quad {\rm as\  elements\
of\ } A. \lb{comp}
\eea
\el

\pr The part $i)$ is a standard consequence of Proposition
\ref{slfrob}.

In order to prove the part $ii)$ use the identity $\phi_L(a_1\ci
(a_2\ast S_A(a_3)))=\phi_L 
((a_1\ast a_3)\ci a_2)$ -- following from (\ref{rigrel}) and the
coherence axioms in $\c$ -- holding true for any $a_1,a_2,a_3\in A$ and
Proposition \ref{slfrob} to show that in ${A^L}\ot {_L A}$ we have
\bea  y_i\ot x_i\ast a &=&
y_i \ot \sla\ci \phi_L ((x_i\ast a)\ci y_j)\ci x_j=
y_i \ci \sla\ci \phi_L ((x_i\ast a)\ci y_j)\ot x_j=\nn
&=&y_i \ci \sla\ci \phi_L (x_i\ci(y_j\ast S_A(a)))\ot x_j=
y_j \ast S_A(a)\ot x_j. \nonumber \eea

The part $iii)$ is checked by direct calculation making use of the
definition (\ref{conva}), the quasi-basis property (\ref{d2}) and the
coherence axioms in $\c$. \hspace{1cm}\qed

We are ready to construct the various ingredients of the Hopf
algebroid structure on $A$. In addition to the map $\sla$ in
(\ref{sl}) -- that is going to be the source map of the left
bialgebroid structure -- with the help of the map $\mu$ in
(\ref{munu}) introduce the map that is going to be the target map
as 
\be  \tla:L^{op} \to \ A \qquad l\mapsto \ \i\x \mu(l). \lb{tl} \ee
It is obviously a ring homomorphism. The relations (\ref{rigrel})
imply  that the element
$ i_A={\coev_L}\ci {\ev_R} \lb{ia} 
\in A$ satisfies
\be \sla\ci \phi_L (a\ci i_A) \ci i_A = a \ci i_A = \tla\ci \phi_L
(a\ci i_A)\ci i_A \lb{lem327}\ee
for any element $a$ in $A$ and for the map $\phi_L$ introduced in (\ref{phil}).

Defining the ring homomorphisms that are to be the source and target
maps of the right bialgebroid structure on $A$ as  
\bea   \sra:&R\to \ A       \qquad& r\mapsto \ \i\x r \lb{sr}\\
       \tra:&R^{op} \to \ A \qquad& r\mapsto \ \nu^{-1}(r)\x \ib \lb{tr}
\eea
and the maps that are going to be the counits as
\bea &\pla:A\to L\qquad &a\mapsto \phi_L(a\ci i_A)\nn
     &\pra:A\to R\qquad &a\mapsto \nu\ci \phi_L(S_A^{-1}(a)\ci i_A)
\lb{plr}\eea
we can prove
\bl \lb{piform}
The maps $\pla$ and $\pra$ together with the maps (\ref{sl}), (\ref{tl})
and (\ref{sr}-\ref{tr}) satisfy
\bea i)&\quad (a\ci i_A)\ast 1_A &=\sla\ci \pla(a)\nn
     ii)&\quad 1_A\ast (a\ci i_A) &=\tla\ci \pla(a)\nn
     iii)&\quad1_A\ast (i_A\ci a) &=\sra\ci \pra(a) \nn
     iv) &\quad(i_A\ci a)\ast 1_A &=\tra\ci \pra(a)
\nonumber\eea
for all $a\in A$.
\el

\pr Using (\ref{phil}) and (\ref{rigrel}) one computes that
\bea \pla(a)&=& \phi_L(a\ci i_A)= 
\r_{\i}\ci (\i\x {\ev_L})\ci \a_{\i,\ib,\i}\ci 
[(a\ci {\coev_L}\ci {\ev_R})\x \i]\ci  \a^{-1}_{\i,\ib,\i}\ci 
(\i \x {\coev_R})\ci \r^{-1}_{\i}=\nn
&=& \r_{\i}\ci (\i\x {\ev_L})\ci \a_{\i,\ib,\i}\ci 
[(a\ci {\coev_L})\x \i]\ci \l_{\i}^{-1}
\lb{pil}\eea

\begin{figure}[h]
\psfrag{a}{\Huge$a$}
\psfrag{lhs}{\Huge$\pla(a)\ =$}
\begin{center}
{\resizebox*{3.5cm}{!}{\includegraphics{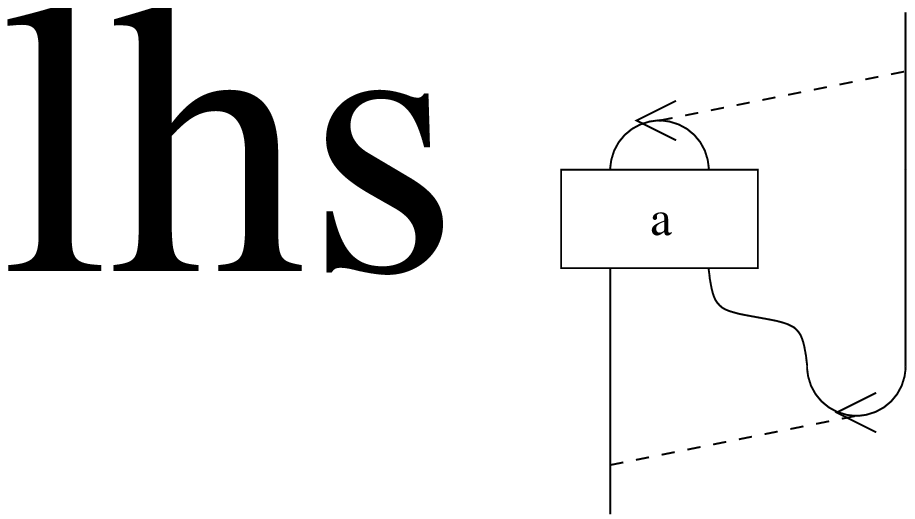}}}
\end{center}
\end{figure}
\vspace{-0.3cm}

\noindent
Then {\em i)} follows by
\bea (a\ci i_A)\ast 1_A &=&
 (\r_{\i}\x \ib)\ci 
[(\i\x {\ev_L})\x \ib]\ci
(\a_{\i,\ib,\i} \x\ib)\ci 
\a^{-1}_{\i\x\ib,\i,\ib} \ci 
[(a\ci {\coev_L}\ci {\ev_R})\x (\i\x\ib)]\ci 
\a_{\i\x\ib,\i,\ib}  \ci \nn
&&(\a^{-1}_{\i,\ib,\i} \x\ib)\ci 
[(\i\x {\coev_R})\x \ib]\ci
 (\r^{-1}_{\i}\x \ib)=\nn
&=&  (\r_{\i}\x \ib)\ci 
[(\i\x {\ev_L})\x \ib]\ci
(\a_{\i,\ib,\i} \x\ib)\ci 
[((a\ci {\coev_L}\ci {\ev_R})\x \i)\x\ib]\ci 
(\a^{-1}_{\i,\ib,\i} \x\ib)\ci \nn
&&[(\i\x {\coev_R})\x \ib]\ci
 (\r^{-1}_{\i}\x \ib)=\nn
&=&  (\r_{\i}\x \ib)\ci 
[(\i\x {\ev_L})\x \ib]\ci
(\a_{\i,\ib,\i} \x\ib)\ci 
[((a\ci {\coev_L})\x \i)\x\ib]\ci 
(\l_{\i}^{-1}\x\ib )=\nn
&=& \pla(a)\x\ib =\sla\ci \pla(a).
\nonumber \eea
Analogously, ii) follows by 
\eject

\begin{figure}[h]
\psfrag{a}{\Huge$a$}
\psfrag{=}{\Huge$=$}
\psfrag{lhs}{\Huge$1_A\ast(a\ci i_A)\ =$}
\psfrag{rhs}{\Huge$=\ \tla\ci \pla (a)$}
\begin{center}
{\resizebox*{4.4cm}{!}{\includegraphics{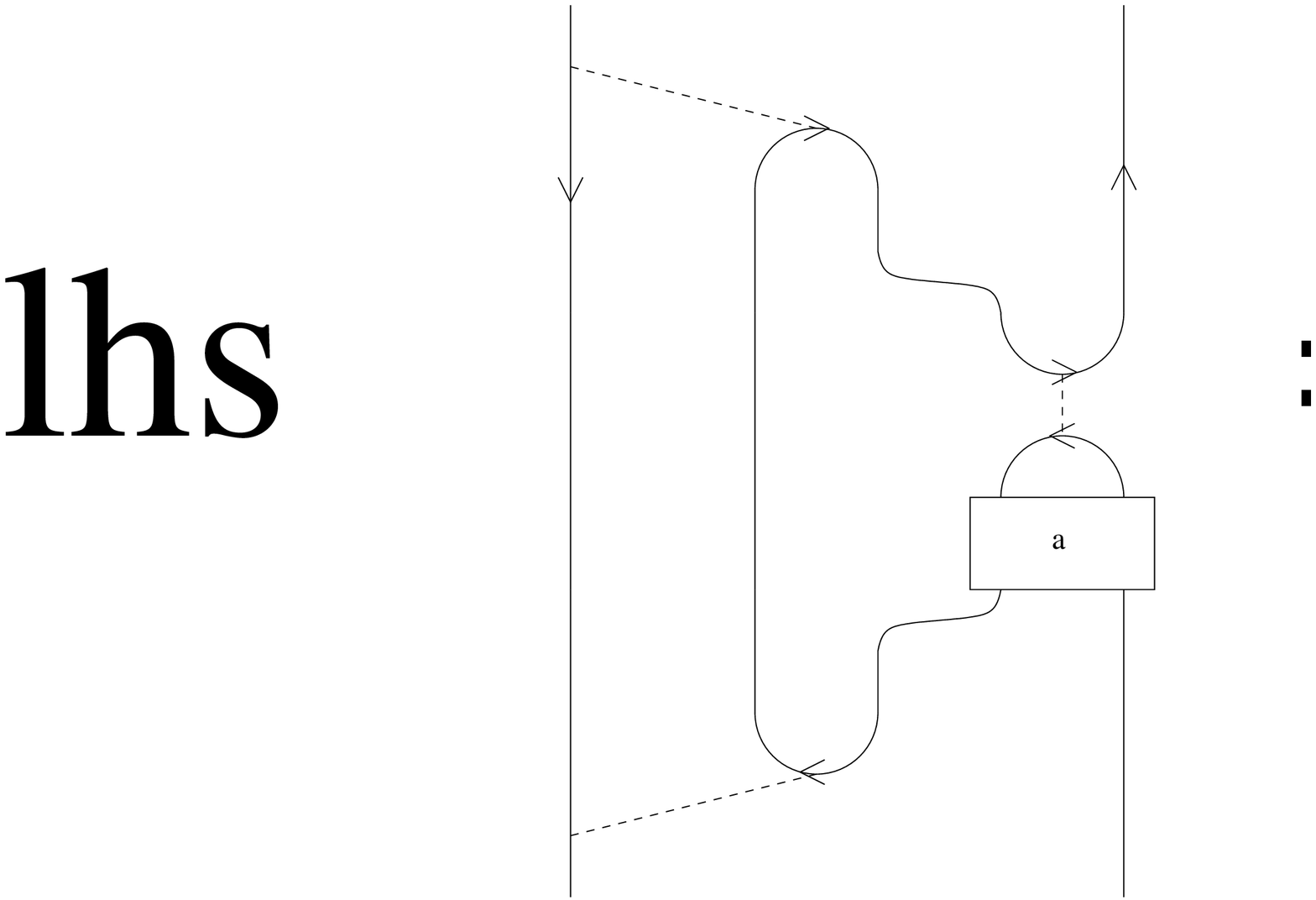}}
~~~~
\resizebox*{4.7cm}{!}{\includegraphics{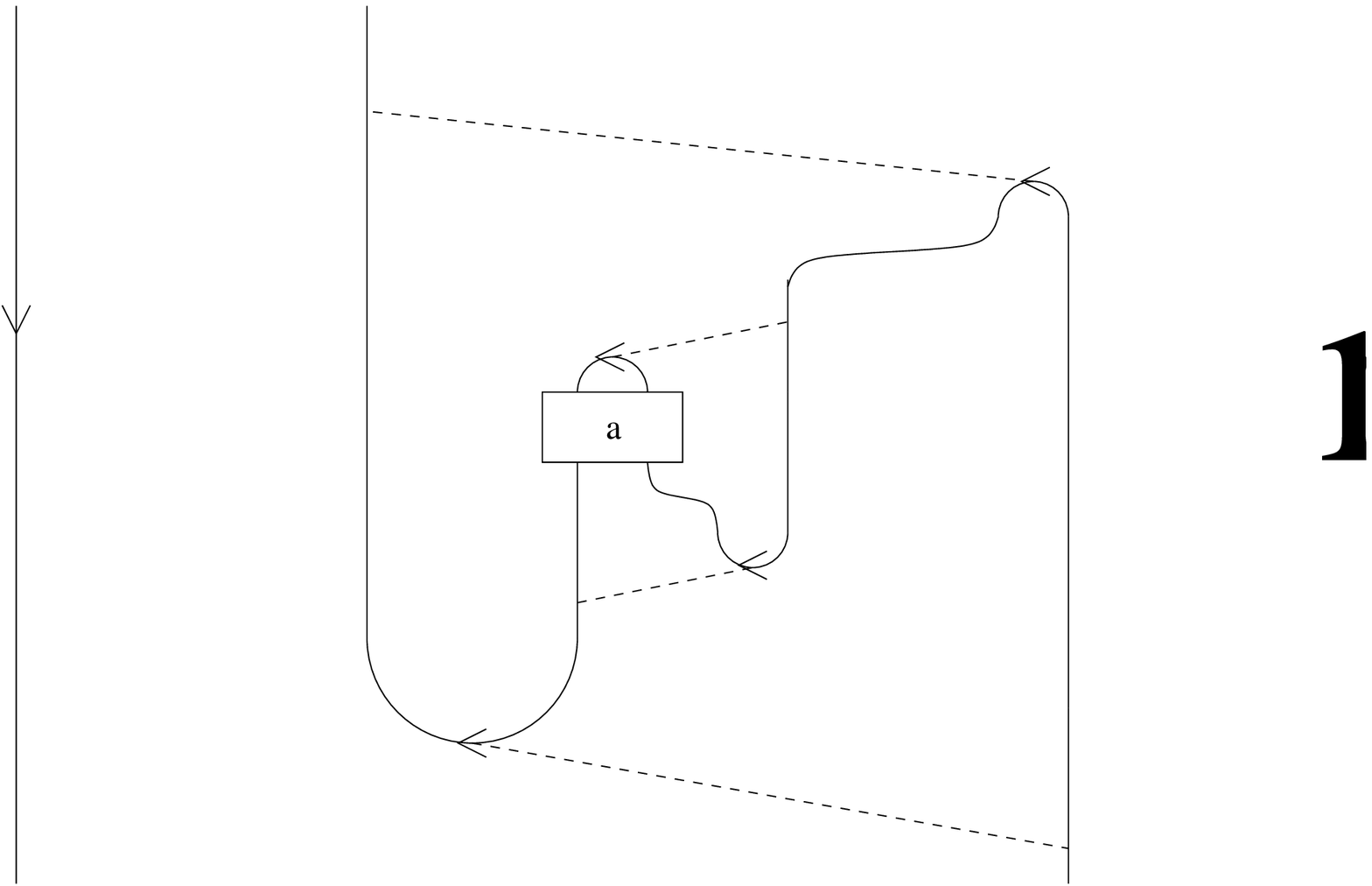}}}
\end{center}
\end{figure}

\noindent
The identities {\em iii)} and {\em iv)} are proven analogously using
\be \pra(a)=\r_{\ib}\ci (\ib \x {\ev_R})\ci (\ib \x a)\ci
\a_{\ib,\i,\ib}\ci ({\coev_R}\x\ib)\ci \l_{\ib}^{-1}
\lb{pir}\ee

\begin{figure}[h]
\psfrag{a}{\Huge$a$}
\psfrag{lhs}{\Huge$\pra(a)\ =$}
\begin{center}
{\resizebox*{3.5cm}{!}{\includegraphics{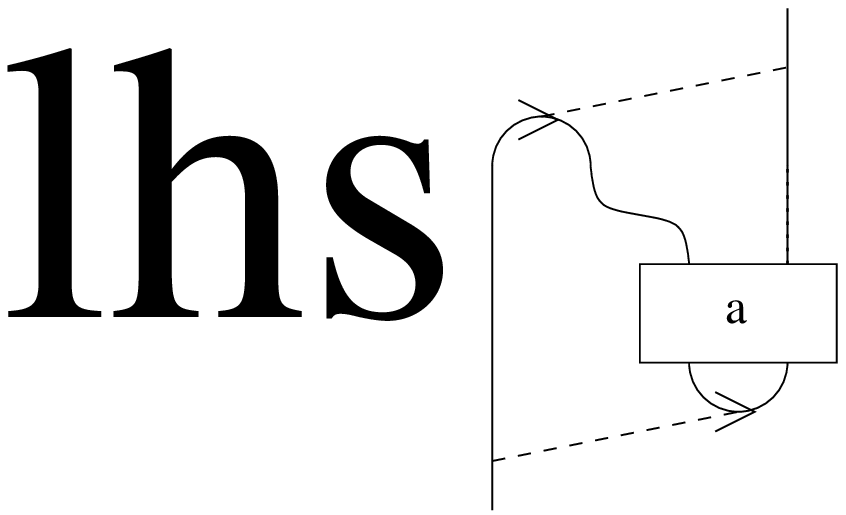}}}
\end{center}
\end{figure}

\qed

\noindent
In the following theorem we give the explicit form of the bialgebroid
structures on the ring $A$ the existence of which follows by Theorem
3.5 in \cite{Sz}:
\bt
The ring $A$ carries a  left bialgebroid structure ${\cal A}_L$  over the
base $L$ and a  right bialgebroid structure  ${\cal A}_R$  over $R$ 
with  the structural maps
\be\begin{array}{lll}
&\sla(l)=l\x\ib \qquad &\tla(l)=\i \x \mu(l) \nn
&\gla(a)=a\ast S_A^{-1}(y_i)\ot x_i = S_A^{-1}(y_i)\ot x_i \ast a\qquad
&\pla(a) = \phi_L(a\ci {i_A})\nn
     &\sra(r)=\i \x r \qquad &\tra(r)= \nu^{-1}(r)\x \ib \nn
&\gra(a)=a\ast S_A(x_i)\ot y_i= S_A(x_i)\ot y_i\ast a \qquad 
&\pra(a)= \nu\ci \phi_L( S_A^{-1}(a)\ci i_A ).
\end{array}\lb{bgd}\ee
\et

\pr The equality of the two forms of $\gla$ given in (\ref{bgd})
follows from (\ref{lem326})  and the fact  that by (\ref{fri}) 
we have $S_A(a_1 \ast a_2)= S_A(a_2)\ast S_A(a_1)$.

Since 
\be (\sla(l)\ci a_1) \ast a_2 = \sla(l) \ci (a_1\ast a_2) \quad {\rm
and}\quad 
     a_1 \ast (\tla(l)\ci a_2) =  \tla(l) \ci (a_1\ast a_2) \lb{*lmod}\ee
the $\gla$ is an $L$-$L$ bimodule map. 
The left $L$-module map property of $\pla$ follows from the left
$L$-module property of $\phi_L$. 
Its right $L$-module map property follows by the use of the identity
${\ev_L}\ci (\mu(l)\x \i)={\ev_L}\ci (\ib\x l)$ and the relations
(\ref{rigrel}): 
\bea \pla(\tla(l)\ci a) &=&  \phi_L(\tla(l)\ci a\ci i_A)=\nn
&=&\r_{\i}\ci (\i\x {\ev_L})\ci \a_{\i,\ib,\i}\ci [(\i \x \ib) \x l]\ci
(a \x \i)\ci  ({\coev_L}\x \i)\ci \l_{\i}^{-1}=\nn
&=&\r_{\i}\ci (\i\x {\ev_L})\ci \a_{\i,\ib,\i}\ci (a \x \i)\ci
({\coev_L}\x \i)\ci \l_{\i}^{-1} \ci l=
\pla(a)\ci l.
\nonumber\eea
The coassociativity of $\gla$ follows by using both forms of it:
\[ (\gla\ot \id_A)\ci \gla(a)= S_A^{-1}(y_j)\ot x_j \ast a
\ast S_A^{-1}(y_i) \ot x_i = (\id_A\ot \gla)\ci \gla(a). \]
The relations (\ref{convmod}-\ref{convmodd}) imply that
$S_A\ci \tla(l)=\sla(l)$.
Since by Proposition \ref{slfrob} and (\ref{lem327}) we have $ y_i \ci
\sla\ci \pla(x_i)=i_A=\sla\ci \pla \ci S_A^{-1} 
(y_i)\ci x_i $ and $i_A$ is the unit for ${\hat A}$, the $\pla$ is the
counit for $\gla$. 

\noindent
Using (\ref{sep}) we have the identity in
$\atl\ot  \asl$ for any $a\in A$ and $l\in L$:
\bea a\di \ci \tla(l) \ot a\dii &=& 
S_A^{-1}(y_i)\ci \tla(l) \ot x_i\ast a =
S_A^{-1}(y_i) \ot (x_i \ci \sla(l))\ast a =\nn
&=& S_A^{-1}(y_i) \ot (x_i \ast a)\ci \sla(l)=
a\di \ot a\dii \ci \sla(l).
\nonumber\eea
Since $a\ast 1_A=\sla\ci \phi_L(a)$ for all $a\in A$, also
\bea \gla(1_A)&=&
S_A^{-1}(y_i)\ot x_i\ast 1_A =
S_A^{-1}(y_i)\ot \sla\ci \phi_L(x_i)=
\tla\ci \phi_L(x_i) \ci S_A^{-1}(y_i)\ot 1_A =
1_A\ot 1_A.
\nonumber\eea
It follows from (\ref{comp}) that
\bea \gla(a_1) \gla(a_2) &=&
S_A^{-1}(y_i)\ci S_A^{-1}(y_j)\ot (x_i\ast a_1)\ci (x_j\ast a_2)=\nn
&=& S_A^{-1}\left( y_k \ci \sla\ci \phi_L (x_k\ci y_j \ci y_i)\right) \ot
(x_i\ast a_1)\ci (x_j\ast a_2)= \nn
&=& S_A^{-1}(y_k)\ot  \sla\ci \phi_L (x_k\ci y_j \ci y_i)\ci  (x_i\ast
a_1)\ci (x_j\ast a_2)= \nn
&=& S_A^{-1}(y_k)\ot [(x_k\ci y_j)\ast a_1]\ci  (x_j\ast a_2)= 
 S_A^{-1}(y_k)\ot x_k\ast (a_1\ci a_2)=
\gla(a_1\ci a_2).
\nonumber\eea
Finally, using (\ref{lem327}) we have
\[ \pla(a_1\ci \sla\ci \pla (a_2)) =
\phi_L(a_1\ci  \sla\ci \phi_L (a_2\ci i_A)\ci {i_A})=
\phi_L(a_1\ci a_2\ci {i_A})=
\pla(a_1\ci a_2)=
\pla(a_1\ci \tla\ci \pla (a_2)).\]
This finishes the proof of the statement that
${\cal A}_L=(A,L,\sla,\tla,\gla,\pla)$ is a left bialgebroid.

In order to prove that ${\cal A}_R=(A,R,\sra,\tra,\gra,\pra)$ is a
right bialgebroid observe that $(S_A,\nu)$ is a left bialgebroid
isomorphism ${\cal A}_L\to ({\cal A}_R)^{op}_{cop}$ i. e.
\[ \begin{array}{ll} 
S_A\ci \sla = \sra\ci \nu \qquad &S_A\ci \tla = \tra \ci \nu \nn
S_{A\ot_L A}\ci \gla = \gra \ci S_A \qquad & \nu \ci
\pla = \pra \ci S_A
\end{array}\]
where the map $S_{A\ot_L A}:\atl\ot \asl\to \asr\ot \atr$ is defined
as $a_1\ot a_2\mapsto S_A(a_2)\ot S_A(a_1)$ and $\nu:L\to R^{op}$ has
been introduced in (\ref{munu}).
\hspace{1cm}\qed


\bt \lb{ahopf}
The left and right bialgebroid structures (\ref{bgd}) on the ring
$A$ and the map $S_A$ in (\ref{s}) form a Hopf algebroid ${\cal A}$.
\et

\pr It is obvious from (\ref{sl}), (\ref{tl}) and (\ref{sr}-\ref{tr})
that  
\[ \sla(L) = \tra(R) \qquad \tla(L)=\sra(R) \]
as subrings of $A$.

Using the explicit forms (\ref{bgd}) of the coproducts $\gla$ and
$\gra$ we have 
\bea (\gla\ot \id_A)\ci \gra(a) &=
S_A^{-1} (y_j)\ot x_j\ast a\ast S_A(x_i) \ot y_i =&
(\id_A\ot \gra)\ci \gla(a)\nn
(\gra\ot \id_A)\ci \gla(a) &=
S_A(x_j)\ot y_j\ast a\ast S_A^{-1}(y_i)\ot x_i =&
(\id_A\ot \gla)\ci \gra(a).
\nonumber \eea
The map $S_A$ is bijective by the bijectivity of the maps $\f$ and
$\fdot$. It is anti-multiplicative by (\ref{fri}) and
(\ref{convmod}-\ref{convmodd}) imply 
$S_A\ci \tla=\sla$ and  $S_A\ci \tra=\sra$. Hence
\[ S_A(\tla(l)\ci  a \ci \tla(l^{\prime}))= \sla(l^{\prime})\ci S_A(a)
\ci \sla(l)
\qquad  S_A(\tra(r)\ci a\ci  \tra(r^{\prime}))= \sra(r^{\prime}) \ci
S_A(a) \ci \sra(r).\] 
By (\ref{comp}), the {\em iii)} and the {\em ii)} of Lemma
\ref{piform} also
\bea S_A(a\di)\ci a\dii &=& y_i\ci (x_i\ast a)= 1_A\ast (i_A\ci a)=
\sra\ci \pra(a)\nn
a\ui\ci S_A(a\uii) &=& S_A \left( (y_i\ast a)\ci x_i \right)= S_A
\left( 1_A\ast (a\ci i_A)\right)= S_A\ci \tla\ci \pla(a)= \sla\ci
\pla(a). \qquad\qed
\nonumber \eea
Theorem \ref{ahopf} above generalizes the result proven in the
Proposition 6.19 in \cite{Mu}. There $\i$ is assumed to be a
finite index and irreducible D2 Frobenius 1-morphism in
a semisimple $k$-linear bicategory. Then the endomorphism rings $A$ of
$\i\x\ib$ and $B$ of $\ib\x \i$ are  equipped with dual finite
dimensional semisimple and cosemisimple Hopf algebra structures.


Interchanging the roles of the 1-morphisms $\i$ and $\ib$, Theorem
\ref{ahopf} implies that also the ring $B$ carries a Hopf
algebroid structure ${\cal B}=({\cal B}_L,{\cal B}_R,S_B)$. The left
bialgebroid ${\cal B}_L$ has $R$ as the base and the structural maps
\be\begin{array}{ll}
\slb(r)=r\x\i \qquad &\tlb(r)=\ib \x \nu^{-1}(r)\nn
\glb(b)=b\ast \f(x_i)\ot \f(y_i)\equiv \f(x_i)\ot \f(y_i)\ast b
\qquad &\plb(b)=\nu\ci \phi_L\ci \f^{-1}(b).\end{array}\lb{bl}\ee
The right  bialgebroid ${\cal B}_R$ has $L$ as the base and the
structural maps 
\be \begin{array}{ll}
\srb(l)=\ib\x l\qquad &\trb(l)\ \ =\ \mu(l)\x\i\nn
\grb(b)=b\ast \fdot(y_i)\ot \fdot(x_i)\equiv  \fdot(y_i)\ot
\fdot(x_i) \ast b\qquad \quad&\prb(b)\ =\ \phi_L\ci \fdot^{-1}(b).
\end{array}\lb{br}\ee
In the rest of this subsection we are going to prove that the Hopf
algebroid ${\cal B}$ is strictly isomorphic to the dual of
${\cal A}$.
\bl The four rings $\asld$ and  $\atld$ in (\ref{ldual}), $\asrd$ and
$\atrd$ in (\ref{rdual}) are all isomorphic to $B$.
\el

\pr We construct the isomorphism

\be  \alpha^*:B\to\asrd\quad\qquad b\mapsto\nu\ci \phi_L(S_A^{-1}(-)\ci
\f^{-1}(b)) \lb{al*}\ee
with inverse 
\be \alpha^{*-1}:\asrd \to B\quad\qquad {\phi^*}\mapsto\f\left(y_i\ci
\tra\ci \phi^*\ci S_A(x_i)\right).\ee
As a matter of fact
\bea \alpha^{*-1}\ci \alpha^* (b)&=&
\f\left( y_i\ci \tra\ci \nu\ci\phi_L(x_i\ci \f^{-1}(b))\right)=b\nn
\alpha^*\ci \alpha^{*-1}(\phi^*)&=& 
\nu\ci \phi_L\left( S_A^{-1}(-)\ci y_i\ci \tra\ci \phi^*\ci
S_A(x_i)\right)=
\nu\ci \phi_L\left( S_A^{-1}(-)\ci y_i\right)\ci \phi^*\ci S_A(x_i)=\nn
&=&\phi^*\ci S_A\left(\sla\ci \phi_L(S_A^{-1}(-)\ci y_i)\ci x_i \right)=
\phi^*\nn
\left({\alpha^*}(b_1){\alpha^*}(b_2)\right)(a)&=&
\nu\ci \phi_L\left(S_A^{-1}[(y_i\ast a)\ci \tra\ci \nu\ci \phi_L(x_i\ci
\f^{-1}(b_2))]\ci \f^{-1}(b_1)\right)=\nn
&=&\nu \ci \phi_L\left( [S_A^{-1}(a)\ast S_A^{-1}\ci \f^{-1}(b_2)]\ci
\f^{-1}(b_1) \right)=\nn
&=&\nu \ci \phi_L\left( S_A^{-1}(a)\ci[\f^{-1}(b_1)\ast \f^{-1}(b_2)]
\right)= 
\nu \ci \phi_L\left( S_A^{-1}(a)\ci \f^{-1}(b_1\ci b_2)]\right)=\nn
&=& {\alpha^*}(b_1\ci b_2)(a).
\nonumber\eea
Analogously, one checks that 
\[\begin{array}{rccccl}
{^*\alpha}:&B&\to&\atrd\quad\qquad b&\mapsto&\nu\ci
\phi_L(\f^{-1}(b)\ci-)\nn
{^*\alpha^{-1}}:&\ \atrd&\to&B\quad\qquad {^*\phi}&\mapsto&\f\left(\tra\ci
{^*\phi}(y_i)\ci x_i\right)\nn
{_*\alpha}:&B&\to& \asld \quad\qquad b&\mapsto& \phi_L(- \ci
\fdot^{-1}(b)\nn
{_*\alpha^{-1}}:& \asld& \to& B \qquad {_*\phi}&\mapsto& 
\fdot\left( y_i\ci \sla \ci {_*\phi}(x_i)\right)\nn
{\alpha_*}:&B& \to& \atld \quad\qquad b&\mapsto& \phi_L(\fdot^{-1}(b)\ci
S_A(-)) \nn
{\alpha_*}^{-1}:&\atld& \to&B\quad\qquad {\phi_*}&\mapsto&
\fdot\left(\sla\ci \phi_*\ci S_A^{-1}(y_i)\ci x_i \right)
\end{array}\]
define isomorphisms of rings.\hspace{1cm}\qed


\bt \lb{int} The element $i_A\colon = {\coev_L}\ci {\ev_R}$ is a two sided 
non-degenerate integral in the Hopf algebroid ${\cal A}$.
\et

\pr The element $i_A$ is a left integral by (\ref{lem327}) and the
definition (\ref{pil}) of the map $\pla$. Since it is invariant under
$S_A$ it is also a right integral. It remains to check
non-degeneracy. As a matter of fact
\[\alpha^*(b) \ru i_A=y_i\ci \tra\ci\nu\ci \phi_L(x_i\ci
\f^{-1}(b))=\f^{-1}(b)\] 
hence the map 
\be {(i_A)_R}= \f^{-1}\ci{\alpha^{*-1}}\lb{als}\ee
is bijective. Analogously, ${_R(i_A)}=\fdot^{-1}\ci {^*\alpha^{-1}}$ is
bijective. 
\hspace{1cm} \qed

Interchanging the roles of the 1-morphisms $\i$ and $\ib$, Theorem
\ref{int} implies that $i_B= {\coev}_R\ci {\ev}_L$ is a two sided
non-degenerate integral in the Hopf algebroid ${\cal B}$.

\br With the help of the two sided non-degenerate integral $i_A=
{\coev_L}\ci {\ev_R}$ the convolution product (\ref{conva}) takes the
forms
\bea a_1\ast a_2 &=& (i_A)_R\left( (i_A)_R^{-1}(a_1) (i_A)_R^{-1}(a_2)
\right)\equiv 
{_R(i_A)}\left( {_R(i_A)}^{-1}(a_2) {_R(i_A)}^{-1}(a_1)
\right)\equiv \nn
&=& (i_A)_L\left( (i_A)_L^{-1}(a_2) (i_A)_L^{-1}(a_1)
\right) \equiv
{_L(i_A)}\left( {_L(i_A)}^{-1}(a_1) {_L(i_A)}^{-1}(a_2)
\right) \eea
where $i_R:\asr\to A$, ${_R i}:\atr\to A$, ${_L i}:\asl \to A$ and 
${i_L}:\atl \to A$ are bijections of additive groups. 
\er

By Theorems \ref{int} and \ref{dualthm}  there exists a Hopf algebroid
structure $\asrd_{i_A}$ on the ring $\asrd$ with structural maps listed in
the Theorem \ref{dualthm}.


\bt The Hopf algebroid  $\asrd_{i_A}$ is strictly isomorphic to ${\cal
B}$. 
\et

\pr  Let $\alpha^*:B\to \asrd$ be the ring isomorphism (\ref{al*}).
We claim that $({\alpha^*},\id_R)$ is a strict isomorphism of Hopf
algebroids ${\cal B}\to \asrd_{i_A}$. Use  
(\ref{convmod}) and (\ref{als}) to check that 
\bea 
{\alpha^{*-1}}\ci s^*_L(r)&=&
\f\left( y_i\ci \tra[r\ci \pra\ci S_A (x_i)]\right)=
\f\left( y_i\ci \sla\ci \pla (x_i) \ci \tra(r) ]\right)=\nn
&=&\f(i_A\ci \tra(r))= r\x\i=\slb(r)\nn
{\alpha^{*-1}}\ci t^*_L(r)&=&
\f \left( y_i\ci \tra\ci \pra(\sra(r)\ci S_A(x_i))\right)=
\f \left( \tra(r)\ci  y_i\ci \sla\ci \pla (x_i)\right)=\nn
&=&\f\left(\tra(r) \ci i_A\right)=\ib\x \nu^{-1}(r)=\tlb(r)\nn
({\alpha^{*-1}}\ot {\alpha^{*-1}})\ci \gamma^*_L\ci
{\alpha^*}(b)&=&
\alpha^{*-1}\left( \alpha^*(b)\lu S_A(x_i)\right)\ot \alpha^{*-1}\ci
(i_A)_R^{-1} (y_i)=\nn
&=&\f\left( x_i\ci y_j\ci \tra\ci \nu\ci \phi_L(x_j\ci
\f^{-1}(b))\right)\ot \f\left(y_i\right)=\nn
&=&\f\left(x_i\ci \f^{-1}(b)\right)\ot \f(y_i)=
\glb(b)\nn 
\pi^*_L\ci \alpha^*(b)&=&
\alpha^*(b)(1_A)=\nu\ci \phi_L\ci \f^{-1}(b)=
\plb(b)\nn
{\alpha^{*-1}}\ci S^*\ci {\alpha^*}&=&
{\alpha^{*-1}}\ci (i_A)_R^{-1}\ci S_A\ci (i_A)_R \ci {\alpha^*}=
\f\ci S_A\ci \f^{-1}= \f\ci \fdot^{-1}=
S_B.\quad \qed\nonumber
\eea

\bex \rm Let $N\to M$ be a D2 Frobenius extension of rings. Then the
$N$-$M$ bimodule ${_N M_M}$ is a D2 Frobenius 1-morphism in the additive
bicategory of bimodules. Applying the above construction we
obtain the Hopf algebroid described in Section 3 of \cite{HGD}.
\eex

\subsection{The inverse construction}
\lb{inv}

In this subsection we address the question what Hopf algebroids arise as
symmetries of abstract D2 Frobenius extensions in the way explaned in
Subsection \ref{depth2}. By Theorem \ref{int} the existence of a two sided
non-degenerate integral is a necessary condition. The main result of
this subsection states that it is also sufficient.

Throughout this subsection let ${\cal H}=({\cal H}_L,{\cal H}_R,S)$ be
a Hopf algebroid with a non-degenerate right integral $\err$. We
use the notation ${\cal H}_L=(H,L,s_L,t_L,\gamma_L,\pi_L)$ and  ${\cal
H}_R=(H,R,s_R,t_R,\gamma_R,\pi_R)$ .
In what follows
we associate a bicategory and a D2 Frobenius 1-morphism of it to
the pair $({\cal H},\err)$. 
The construction is built on the generalization of the result in
\cite{Y} described in the Appendix. We arrive to the statement that if
$\err$ is a {\em two sided} non-degenerate integral then the
Hopf algebroid symmetry of the D2 Frobenius 1-morphism constructed is
isomorphic to ${\cal H}$.

\smallskip


Recall that for a right bialgebroid ${\cal H}_R$ the right regular
$H$-module $H_H$ is the object part of a comonoid. That is the triple
$(H_H,\gamma_R,\pi_R)$ is a comonoid in the monoidal category ${\cal
M}_H$. Now we claim that for a Hopf algebroid ${\cal H}$ possessing a
non-degenerate right integral the $H_H$ has more structure then just being
a comonoid. The following proposition generalizes the result of
\cite{Mu} on Hopf algebras:

\bp Let ${\cal H}$ be  a Hopf algebroid with a non-derenerate right
integral $\err$. Define the convolution product on the additive
group ${\tt H}$ underlying the ring $H$ as
\be \lb{convh} h_1\ast h_2 =  {_L\err} ({_L\err}^{-1}(h_1)
{_L\err}^{-1}(h_2) )\ee 
and the map 
\be \lb{convun}\eta: R\to H \qquad r\mapsto \err s_R(r) \equiv  \err
t_R(r).\ee 
Then $(H_H, \ast,\eta,\gamma_R,\pi_R)$ is a Frobenius algebra in the
category ${\cal M}_H$. 
\ep

\pr Let ${_*\rho}\colon = {_L\err}^{-1}(1_H)$.  We use the notation
$\gamma_L(h)=h\di\ot h\dii$ 
and $\gamma_R(h)=h\ui\ot h\uii$ for $h\in H$. Using the  Hopf
algebroid identity  $h\ld{_L i}^{-1}(k)=k\lu {i_L}^{-1}(h)$ for all
$h,k\in H$ -- see \cite{HGD} Lemma 5.11 -- the convolution product
(\ref{convh}) has the equivalent forms
\bea h\ast k &=& 
h\ld {_L\err}^{-1}(k)\equiv 
t_L\ci {_*\rho}\left(h\dii S(k)\right) h\di=\nn
&=&k\lu {\err_L }^{-1} (h)\equiv
s_L\ci {_*\rho}\left( h S(k\di)\right) k\dii.
\lb{convform}\eea
First we show that $(H_H,\ast,\eta)$ is a monoid in ${\cal M}_H$.
Both $\ast$ and $\eta$ are $H$-module
maps as for all $h,k,k^{\prime}\in H$ and $r\in R$ we have
\bea k h\ui\ast k^{\prime}h\uii &=& 
t_L\ci {_*\rho}\left( k\dii {h\ui}\dii S( k^{\prime} h\uii)\right) 
k\di {h\ui}\di=
 t_L\ci {_*\rho}\left( k\dii s_L\ci \pi_L( h\dii) S(k^{\prime})
\right) 
k\di h\di=\nn
&=& t_L\ci {_*\rho}\left( k\dii S(k^{\prime}) \right) k\di h=
(k\ast k^{\prime})h\nn
\eta\ci \pi_R\left (s_R(r) h\right)&=&
\err s_R\ci \pi_R \left (s_R(r) h\right) = 
\err s_R(r) h=\eta(r)h.
\nonumber\eea
Using (\ref{convform}), the convolution product $\ast$ can be seen to
be associative: for $h,k,l\in H$ we have
\[ (h\ast k)\ast l =
t_L\ci {_*\rho} \left( k_{(3)} S(l)\right) s_L\ci {_*\rho}\left( h
S(k\di)\right) k\dii=
h\ast ( k\ast l). \]
As $i$ is the unit for the ring $({\tt H},\ast)$, the map $\eta$ is
the unit for $\ast$: 
\bea h\ast \eta(r)=& h\ast (\err s_R(r))=(h\ast
\err)s_R(r)=hs_R(r)&=r_H (h\ot r)\nn
\eta(r)\ast h =& (\err t_R(r))\ast h = (\err \ast h)t_R(r)=h t_R(r)&=
l_H (r\ot h).
\nonumber\eea
The compatibility of the monoid $(H_H,\ast,\eta)$ and the comonoid
$(H_H,\gamma_R,\pi_R)$ follows by
\bea \gamma_R(h\ast k) &=& 
\gamma_R(h \ld {_L\err}^{-1}(k))=
h\ui\ot h\uii \ld {_L\err}^{-1}(k) =
h\ui\ot h\uii \ast k \nn
&=& \gamma_R(k\lu {\err_L}^{-1}(h))=
k\ui \lu {\err_L}^{-1}(h)\ot k\uii=
h\ast k\ui \ot k\uii . \qquad \qed
\nonumber\eea
It is obvious from (\ref{convform}) that also
\be \gamma_L(h\ast k)=h\di \ot h\dii\ast k=h\ast k\di \ot
k\dii. \ee


\bigskip

Let $\BIM({\cal M}_H)$ denote the  additive bicategory  of
internal bimodules in  the monoidal category ${\cal M}_H$
of right $H$-modules. (For the definition of the bicategory  of
internal bimodules see the Appendix.) Denote the trivial monoid
corresponding to the monoidal unit $R_H$ of ${\cal M}_H$ by  ${\cal
U}=(R_H,l_R=r_R,R)$ and the monoid (\ref{convh}-\ref{convun}) by
${\cal Q}=(H_H,\ast,\eta)$. Then by Proposition \ref{app} the 
${\cal X}\colon = (H_H,l_H,\ast)$  is an
internal  ${\cal U}$-${\cal Q}$ bimodule and  $\overline{{\cal X}}\colon
= (H_H,\ast,r_H)$ is an internal  
${\cal Q}$-${\cal U}$ bimodule in ${\cal M}_H$.
Furthermore, $\overline{{\cal X}}$ is the two sided dual of ${\cal X}$
in  $\BIM({\cal M}_H)$.

\bp \lb{ad2}
The Frobenius 1-morphism 
$\ix$ of the  bicategory $\BIM({\cal M}_H)$ satisfies the D2
condition. 
\ep

\pr The 1-morphism $\ix \stac{{\cal Q}} \ixb$ of $\BIM({\cal
M}_H)$ is the internal ${\cal U}$-${\cal U}$ bimodule $(H_H,l_H,r_H)$.
Introduce the left multiplication map
\be\lb{Lambda} \Lambda: H\to {\cal M}_H^1(H_H,H_H)
\equiv \BIM({\cal M}_H)^2(
\ix \stac{{\cal Q}} \ixb,\ix \stac{{\cal Q}} \ixb)
\qquad  \Lambda(h)k\colon = hk.\ee
We construct the D2 quasi-basis 
\be \Lambda\left(S(i\di)\right)\ot 
    \Lambda\left(i\dii \right) .\ee
Using the explicit forms of the coherence isomorphisms
\be \begin{array}{rlrl}
     \underline{\l}_{\ix}(r\ot h)&=\ ht_R(r)&\qquad
     \underline{\l}_{\ixb}(h\ot k)\ =& h\ast k\nn
     \underline{\r}_{\ix}(h\ot k)&=\ h\ast k&\qquad
     \underline{\r}_{\ixb}(h\ot r)\ =& h s_R(r)\nn
     \underline{\a}_{\ix,\ixb,\ix}(h\ot k\ot l)&=\ h\ot k\ot l&
\qquad\qquad
     \underline{\a}_{\ixb,\ix,\ixb}(h\ot k\ot l)\ =&h\ot k\ot l
\end{array}\lb{coh}\ee
for $r\in R$, $h,k,l\in H$ and the 2-morphisms (\ref{left}) and
(\ref{right})  one checks that
\bea (\Lambda\left(S(i\di)\right)\stac{{\cal U}} \ix)&\ci&
\gamma_R\ \ci\  \ast\ \ci \ 
(\Lambda\left(i\dii \right)\stac{{\cal U}} \ix)
(h\ot k)=\nn
&=&S(i\di)[(i\dii h)\ast k]\ui \ot [(i\dii h)\ast k]\uii=\nn
&=&h  S(i\di) {i\dii}\ui \ot {i\dii}\uii\ast k=
h s_R\ci \pi_R (i\ui) \ot i\uii\ast k=\nn
&=&h\ot ( i\uii\ast k)t_R\ci \pi_R (i\ui)=
h\ot k \nonumber
\eea
which is the D2 quasi-basis property (\ref{d2}). \hspace{1cm} \qed

Proposition \ref{ad2} together with Theorem \ref{ahopf} implies that
the ring 
$ A\colon =\BIM({\cal M}_H )^2(\ix \stac{{\cal Q}} \ixb,\ix \stac{{\cal
Q}} \ixb)$ carries a Hopf algebriod structure ${\cal A}$.


\bt If the non-degenerate right integral $\err$ used to define the
convolution product (\ref{convh}) is invariant under the antipode
then the Hopf algebroid ${\cal A}$
is  isomorphic to ${\cal H}$.  
\et

\br The non-degenerate {\em right} integral is also a {\em left}
integral if and only if it is invariant under $S_A$. The {\em if} part
follows from {\em iii)} of Lemma \ref{lemint}, and the {\em only if}
part follows by
\[S_A^{-1}(i)=i\ld \left( i \rd {_*\rho}\right)=
t_L\ci  {_*\rho}(i\dii i)i\di=
t_L\ci  {_*\rho}(s_L\ci \pi_L(i\dii) i)i\di=
t_L\ci  {_*\rho}(i)t_L\ci \pi_L(i\dii)i\di=i.
\]
\er

{\em Proof of the Theorem:} The total ring of the Hopf algebroid
${\cal A}$ is
\[ A\colon =
\BIM({\cal M}_H )^2(\ix \stac{{\cal Q}} \ixb,\ix \stac{{\cal Q}} \ixb)
\equiv {\cal M}_H^1(H_H,H_H)\simeq H \]
via the ring isomorphism defined by the left
multiplication map $\Lambda:H\to A$ as in (\ref{Lambda}).
 
The base ring of the left bialgebroid underlying ${\cal A}$ is
$L^A =  \BIM({\cal M}_H )^2(\ix,\ix)$ that is
\[ L^A=\{ \kappa\in {\cal M}_H^1 (H_H,H_H) \ \vert\ 
\kappa(h_1\ast h_2)=
\kappa(h_1)\ast h_2\quad \forall h_1,h_2\in H\ \} \simeq L \]
via the ring isomorphism $\Lambda\ci s_L: L\to L^A$.

We claim that $(\Lambda: H\to A,\Lambda\ci s_L:L\to L^A)$ is a strict isomorphism
of Hopf algebroids. Denote  the
convolution product (\ref{convh}) by $\ast$  and the one in
(\ref{conva}) by $\ast_A$. Also write $S$ for the antipode in ${\cal
H}$ and $S_A$ for the one in ${\cal A}$.
Substituting the 2-morphisms (\ref{coh}), (\ref{left}) 
(\ref{right}) we obtain for all $h,k,m\in H$ the identities
\bea S_A \ci \Lambda(h) (m)&=& 
i\ui s_R\ci \pi_R\left((hi\uii)\ast m\right)=
\Lambda(S(h))(m) \nn
\left[\Lambda(h)\ast_A\Lambda(k)\right](m)&=&
hm\ui \ast km\uii=
\Lambda(h\ast k) (m), \lb{Phiconv}\eea
where in the second step of the first line we used the {\em left}
integral property of $i$. 
The identities (\ref{Phiconv}) imply 
\[\gla\ci \Lambda(h)= 
\Lambda(h)\ast_A S_A^{-1}\ci \Lambda (S(i\di))\ot \Lambda(i\dii)=
\Lambda(h\ast i\di)\ot \Lambda(i\dii)=
(\Lambda\ot\Lambda)\ci \gamma_L(h).\]
Computing the value of $\mu_A$ introduced in (\ref{munu}) on the
element $\Lambda\ci s_L(l)$ of $L^A$, we have
$\mu_A(\Lambda\ci s_L(l))(h)=t_L(l)h$ for all $h\in H$ ,
hence
\bea \sla\ci \Lambda\ci s_L(l)&=&\Lambda\ci s_L(l) \nn
\tla\ci \Lambda\ci s_L(l)&=& \mu_A\ci \Lambda\ci s_L(l)=
\Lambda\ci t_L(l).\nonumber\eea
Finally, for $h,k\in H$
\[ \left[\pla\ci \Lambda(h)\right](k)=
hi\ast k 
=\left[\Lambda\ci s_L\ci \pi_L(h)\right](k).\qquad\qed\]

\vfill
\eject

\section{Appendix: The bicategory of internal bimodules}


Let $({\cal M},\ot,U,\l,\r,\a)$ be a monoidal category with
coequalizers s.t. the monoidal product preserves the
coequalizers. Then there is a bicategory $\BIM({\cal M})$
of internal bimodules constructed as follows.

The 0-morphisms are the monoids in ${\cal M}$, the 1-morphisms with
source $N$ and target $M$ the $M$-$N$ bimodules in ${\cal M}$ and the
2-morphisms the bimodule maps in ${\cal M}$. The vertical composition is
given by the composition in ${\cal M}$. The horizontal composition is the
tensor product of bimodules defined with the help of a coequalizer.

Let ${\cal R}=(R,m_R,\eta_R)$, ${\cal S}=(S,m_S,\eta_S)$ and  ${\cal
T}=(T,m_T,\eta_T)$ be monoids, $(M,\lambda_M, \rho_M)$ an ${\cal R}$-${\cal S}$ and
$(N,\lambda_N,\rho_N)$ an ${\cal S}$-${\cal T}$-bimodule in ${\cal
M}$. As an object in ${\cal M}$ the $M\stackrel{\ot}{_{_{\cal S}}} N$ be
the object part of the coequalizer $(\tau_{M,N}, M\stackrel{\ot}{_{_{\cal S}}} N)$:

\begin{figure}[h]
\psfrag{p1}{$M\ot(S\ot N)$}
\psfrag{p2}{$M\ot N$}
\psfrag{p3}{$M \stackrel{\ot}{_{_{\cal S}}} N$}
\psfrag{a1u}{$(\rho_M\ot N)\ci \a^{-1}_{M,S,N}$}
\psfrag{a1d}{$M\ot\lambda_N$}
\psfrag{a2}{$\tau_{M,N}$}
\begin{center}
{\resizebox*{11cm}{!}{\includegraphics{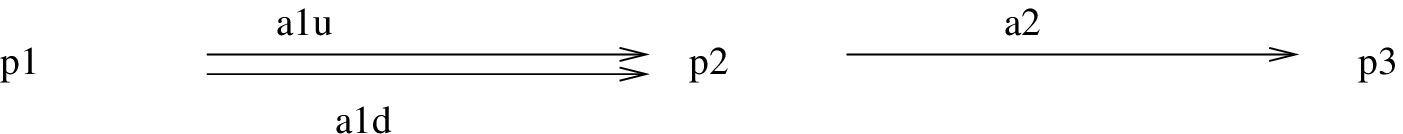}}}
\end{center}
\end{figure}

\noindent

The object $ M\stackrel{\ot}{_{_{\cal S}}} N$  can be equipped
with an ${\cal R}$-${\cal T}$ bimodule structure using the universal
property of the coequalizer:

\begin{figure}[h]
\psfrag{p1}{$R\ot(M\ot(S\ot N))$}
\psfrag{p2}{$R\ot(M\ot N)$}
\psfrag{p3}{$R\ot(M\stackrel{\ot}{_{_{\cal S}}} N)$}
\psfrag{p4}{$M\stackrel{\ot}{_{_{\cal S}}} N$}
\psfrag{a1u}{$R\ot (\rho_M\ot N)\ci \a^{-1}_{M,S,N}$}
\psfrag{a1d}{$R\ot(M\ot\lambda_N)$}
\psfrag{a2}{$R\ot \tau_{M,N}$}
\psfrag{a3}{$\tau_{M,N}\ci (\lambda_M\ot N)\ci\a^{-1}_{R,M,N}$}
\psfrag{a4}{$\lambda_{M\stackrel{\ot}{_{_{\cal S}}} N}$}
\begin{center}
{\resizebox*{11cm}{!}{\includegraphics{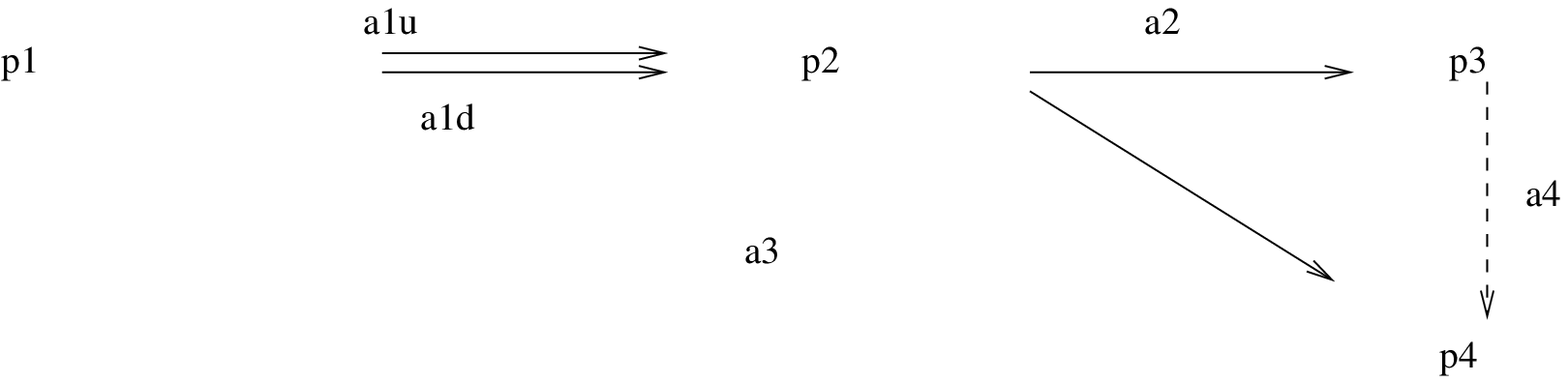}}}
\end{center}

\psfrag{p1}{$(M\ot(S\ot N))\ot T$}
\psfrag{p2}{$(M\ot N)\ot T$}
\psfrag{p3}{$(M\stackrel{\ot}{_{_{\cal S}}} N)\ot T$}
\psfrag{p4}{$M\stackrel{\ot}{_{_{\cal S}}} N$}
\psfrag{a1u}{$(\rho_M\ot N)\ci \a^{-1}_{M,S,N}\ot T$}
\psfrag{a1d}{$(M\ot\lambda_N)\ot T$}
\psfrag{a2}{$\tau_{M,N}\ot T$}
\psfrag{a3}{$\tau_{M,N}\ci (M\ot\rho_N)\ci\a_{M,N,T}$}
\psfrag{a4}{$\rho_{M\stackrel{\ot}{_{_{\cal S}}} N}$}
\begin{center}
{\resizebox*{11cm}{!}{\includegraphics{lambda.eps}}}
\end{center}
\end{figure}

\noindent
One checks that $(M\stackrel{\ot}{_{_{\cal S}}}
N,\lambda_{M\stackrel{\ot}{_{_{\cal S}}} N},
\rho_{M\stackrel{\ot}{_{_{\cal S}}} N} )$ is an ${\cal R}$-${\cal T}$
bimodule in ${\cal M}$.

For the 2-morphisms $p:M\to M^{\prime}$ and $q:N\to
N^{\prime} $ the ${\cal S}$-module product
$p \stackrel{\ot}{_{_{\cal S}}} q$  is constructed also using
the universality of the coequalizer:

\begin{figure}[h]
\psfrag{p1}{$M\ot(S\ot N)$}
\psfrag{p2}{$M\ot N$}
\psfrag{p3}{$M\stackrel{\ot}{_{_{\cal S}}} N$}
\psfrag{p4}{$M^{\prime}\stackrel{\ot}{_{_{\cal S}}} N^{\prime}$}
\psfrag{a1u}{$(\rho_M\ot N)\ci \a^{-1}_{M,S,N}$}
\psfrag{a1d}{$M\ot\lambda_N$}
\psfrag{a2}{$\tau_{M,N}$}
\psfrag{a3}{$\qquad\qquad\qquad\tau_{M^{\prime},N^{\prime}}\ci (p
\ot q)$} 
\psfrag{a4}{$p\stackrel{\ot}{_{_{\cal S}}} q$}
\begin{center}
{\resizebox*{11cm}{!}{\includegraphics{lambda.eps}}}
\end{center}
\end{figure}

\noindent
It is straightforward to check that $\stackrel{\ot}{_{_{\cal S}}}$ is
functorial in both arguments and that the chosen coequalizer $\tau_{M,N}$
becomes a natural transformation from $\ot$ to $\stackrel{\ot}{_{_S}}$.

The coherence natural isomorphisms $\underline{\a}$, $\underline{\l}$ and
$\underline{\r}$ are constructed also using the universality of the
coequalizer:
\vfill
\eject

\begin{figure}[h]
\psfrag{p1}{$M\ot(S\ot S)$}
\psfrag{p2}{$M\ot S$}
\psfrag{p3}{$M\stackrel{\ot}{_{_{\cal S}}} S$}
\psfrag{p4}{$\quad M$}
\psfrag{a1u}{$(\rho_M\ot S)\ci \a^{-1}_{M,S,S}$}
\psfrag{a1d}{$M\ot m_S$}
\psfrag{a2}{$\tau_{M,S}$}
\psfrag{a3}{$\qquad\qquad\qquad \qquad\qquad\qquad\rho_M$} 
\psfrag{a4}{$\underline{\r}_M$}
\begin{center}
{\resizebox*{11cm}{!}{\includegraphics{lambda.eps}}}
\end{center}

\psfrag{p1}{$S\ot(S\ot N)$}
\psfrag{p2}{$S\ot N$}
\psfrag{p3}{$S\stackrel{\ot}{_{_{\cal S}}} N$}
\psfrag{p4}{$\quad N$}
\psfrag{a1u}{$(m_S\ot N)\ci \a^{-1}_{S,S,N}$}
\psfrag{a1d}{$S\ot \lambda_N$}
\psfrag{a2}{$\tau_{S,N}$}
\psfrag{a3}{$\qquad\qquad\qquad \qquad\qquad\qquad\lambda_N$} 
\psfrag{a4}{$\underline{\l}_N$}
\begin{center}
{\resizebox*{11cm}{!}{\includegraphics{lambda.eps}}}
\end{center}
\end{figure}

\begin{figure}[h]
\psfrag{p1}{$(((M\ot S)\ot N)\ot T)\ot Q$}
\psfrag{p2}{$(M\ot N)\ot Q$}
\psfrag{p3}{$(M\stackrel{\ot}{_{_{\cal S}}} N)\stackrel{\ot}{_{_{\cal
T}}} Q$}
\psfrag{p4}{$M\stackrel{\ot}{_{_{\cal S}}}
(N\stackrel{\ot}{_{_{\cal T}}} Q)$}
\psfrag{a1u}{$((M\ot \rho_N)\ot Q)\ci (\a_{M,N,T}\ot Q)
\ci (((\rho_M\ot N)\ot T)\ot Q)$}
\psfrag{a1d}{$((M\ot N) \ot \lambda_Q)\ci \a_{M\ot N,T,Q}\ci
(((M\ot \lambda_N)\ot T)\ot Q)\ci ((\a_{M,N,S}\ot T)\ot Q)$}
\psfrag{a2}{$\tau_{M\stac{{\cal S}} N,Q}\ci (\tau_{M,N}\ot Q)$}
\psfrag{a3}{$\tau_{M,N\stackrel{\ot}{_{_{\cal T}}} Q}\ci (M\ot
\tau_{N,Q}))\ci \a_{M,N,Q}$} 
\psfrag{a4}{$\underline{\a}_{M,N,Q}$}
\begin{center}
{\resizebox*{15cm}{!}{\includegraphics{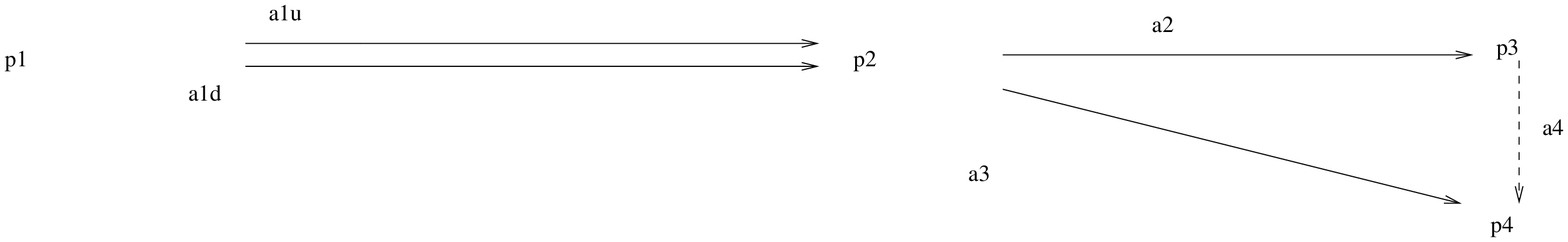}}}
\end{center}
\end{figure}
\smallskip

\noindent
One checks that all $\underline{\r}_M$, $\underline{\l}_N$ and
$\underline{\a}_{M,N,Q}$ are iso's, they are natural and satisfy
the coherence axioms.
This finishes the construction of the bicategory $\BIM({\cal
M})$. If the monoidal category ${\cal M}$ is additive then so is
the bicategory  $\BIM({\cal M})$.


\bp \lb{app} Let  $({\cal M},\ot,U,\l,\r,\a)$ be a monoidal category with
coequalizers s.t. the monoidal product preserves the
coequalizers. Then
\begin{trivlist}
\item[i)] ${\cal U}=(U,\l_U=\r_U,U)$ is a monoid in ${\cal M}$
\item[ii)] For any object $X$ of ${\cal M}$ the ${_{\cal U} X_{\cal
U}}=(X,\l_X,\r_X)$ is an 
${\cal U}$-${\cal U}$ bimodule in ${\cal M}$. Any morphism of ${\cal M}$
is an ${\cal U}$-${\cal U}$ bimodule map. The tensor product over ${\cal U}$ 
in $\BIM({\cal M})$ coincides with the monoidal product
of ${\cal M}$.
\item[iii)] For a monoid ${\cal Q}=(Q,m,\eta)$ in ${\cal M}$ the
\hspace{1cm}\begin{trivlist}
\item[] \hspace{1cm}a) ${_{\cal Q} Q_{\cal Q}}=(Q,m,m )$ is a ${\cal Q}$-${\cal Q}$ bimodule
in ${\cal M}$
\item[] \hspace{1cm}b) ${_{\cal U} Q_{\cal Q}}=(Q,\l_Q,m)$ is a ${\cal U}$-${\cal Q}$ bimodule
in ${\cal M}$
\item[] \hspace{1cm}c) ${_{\cal Q} Q_{\cal U}}=(Q,m,\r_Q)$ is a ${\cal Q}$-${\cal U}$ bimodule
in ${\cal M}$
\item[] \hspace{1cm}d) The 1-morphism  ${_{\cal Q} Q_{\cal U}}$ is the left dual of
${_{\cal U} Q_{\cal Q}}$ in the bicategory $\BIM({\cal M})$
\end{trivlist}
\item[iv)] For a Frobenius algebra $({\cal Q},\delta,\epsilon)$ in ${\cal M}$
the 1-morphism  ${_{\cal Q} Q_{\cal U}}$ is the two sided  dual of
${_{\cal U} Q_{\cal Q}}$ in the bicategory $\BIM({\cal M})$.
\end{trivlist}
\ep

The part iv) is a generalization of Yamagami's analogous result
\cite{Y} where the same claim is proven under the additional
assumption of split separability of the Frobenius algebra  $({\cal
Q},\delta,\epsilon)$.

\pr The only non-trivial part of the statement is {\em d)} of {\em
iii)}  and {\em iv)}.
In order to prove {\em d)} of {\em iii)} we construct the 2-morphisms
\bea \BIM({\cal M})^2({_{\cal Q} Q_{\cal U}}\stac{{\cal U}}
{_{\cal U} Q_{\cal Q}} , {_{\cal Q} Q_{\cal Q}}) &\ni\ \quad {\ev_L}&\colon =
\ m \nn
\BIM({\cal M})^2({_{\cal U} U_{\cal U}}, {_{\cal U} Q_{\cal Q}} 
\stac{{\cal Q}} {_{\cal Q} Q_{\cal U}}) &\ni\ {\coev_L}&\colon = \ 
t^{-1}\ci \eta
\lb{left}\eea
where $t$ is the 2-morphism $\underline{\r}_{_{\cal Q} Q_{\cal
Q}}=\underline{\l}_{_{\cal Q} Q_{\cal Q}}$ regarded as a 2-morphism in
$\BIM({\cal M})( {_{\cal U} Q_{\cal Q}} \stac{{\cal Q}} {_{\cal
Q} Q_{\cal U}}, {_{\cal U} Q_{\cal U}})$.
The ${\ev_L}=m$ is a ${\cal Q}$-${\cal Q}$ bimodule map by the
associativity of $m$.

\noindent
By the unit  property of $\eta$ the first two relations in
(\ref{rigrel}) hold true. 

Similarly, {\em iv)} follows by the existence of 2-morphisms
\bea \BIM({\cal M})^2({_{\cal U} Q_{\cal Q}} 
\stac{{\cal Q}} {_{\cal Q} Q_{\cal U}},{_{\cal U} U_{\cal U}})  &\ni\
\quad {\ev_R}&\colon = \  \varepsilon \ci t
\nn
\BIM({\cal M})^2( {_{\cal Q} Q_{\cal Q}},{_{\cal Q} Q_{\cal
U}}\stac{{\cal U}} {_{\cal U} Q_{\cal Q}})  &\ni\
{\coev_R}&\colon = \  \delta.
\lb{right}\eea
The ${\coev_R}=\delta$ is a ${\cal Q}$-${\cal Q}$ bimodule map by the
Frobenius algebra property. By the counit property of $\varepsilon$ the last two
relations in (\ref{rigrel}) hold true. 
\hspace{1cm}\qed

\end{document}